\input ./style/arxiv-general.cfg
\documentclass[aos,MSNbibl,nameyear,dvips]{arximspdf}
\makeatletter
   \@ifpackageloaded{graphicx}{}{\usepackage{graphicx}}
\makeatother
\usepackage{accents}

\doi{10.1214/14-AOS1302}
\volume{43}
\issue{3}
\pubyear{2015}
\firstpage{1141}
\lastpage{1177}
\docsubty{FLA}

\makeatletter
\renewcommand{\l}{\ell}
\newcommand{\telque}{:}
\newcommand{\Fbar}{ \ol{F}}
\newcommand{\pnorm}{ \ol{\Phi}}

\newcommand{\eqref}[1]{(\ref{#1})}
\newtheorem{theorem}{Theorem}[section]
\newtheorem{proposition}[theorem]{Proposition}

\newtheorem{lemma}[theorem]{Lemma}

\newproclaim{remark}[theorem]{Remark}
\newproclaim{example}[theorem]{Example}
\newproclaim{definition}[theorem]{Definition}
\newproclaim{fact}[theorem]{Fact}

\newcommand{\R}{\mathbb{R}} 
\renewcommand{\P}{\mathbb{P}}
\newcommand{\E}{\mathbb{E}} 
\newcommand{\cond}{ | }
\newcommand{\mtc}{\mathcal}
\newcommand{\mbf}{\mathbf}
\newcommand{\wt}[1]{{\tilde{#1}}}
\newcommand{\wh}[1]{{\hat{#1}}}
\newcommand{\ol}[1]{\bar{#1}}
\newcommand{\FDP}{\operatorname{\textsc{FDP}}}
\makeatother

\begin{document}
\begin{frontmatter}

\title{New procedures controlling the false discovery proportion
via Romano--Wolf's heuristic}
\runtitle{New procedures controlling the FDP}

\begin{aug}
\author[A]{\fnms{Sylvain}~\snm{Delattre}\ead
[label=e1]{sylvain.delattre@univ-paris-diderot.fr}}
\and
\author[B]{\fnms{Etienne}~\snm{Roquain}\corref{}\ead
[label=e2]{etienne.roquain@upmc.fr}\thanksref{T1}}
\runauthor{S. Delattre and E. Roquain}
\thankstext{T1}{Supported by the French Agence Nationale de la
Recherche (ANR Grant references: ANR-09-JCJC-0027-01, ANR-PARCIMONIE,
ANR-09-JCJC-0101-01).}
\affiliation{Universit\'{e} Paris Diderot and Sorbonne Universit\'{e}s}
\address[A]{Univ. Paris Diderot\\
Sorbonne Paris Cit\'{e}\\
LPMA, UMR 7599, F-75205, Paris\\
France\\
\printead{e1}}
\address[B]{Sorbonne Universit\'{e}s\\
UPMC Univ. Paris 6\\
UMR 7599, LPMA, F-75005, Paris\\
France\\
\printead{e2}}

\end{aug}

%
\received{\smonth{6} \syear{2014}}
%
\revised{\smonth{12} \syear{2014}}

%
\begin{abstract}
The false discovery proportion (FDP) is a convenient way to account for
false positives when a large number $m$ of tests are performed simultaneously.
Romano and Wolf [\emph{Ann. Statist.} \textbf{35} (2007) 1378--1408]
have proposed a general principle that builds FDP controlling
procedures from $k$-family-wise error rate controlling procedures while
incorporating dependencies in an appropriate manner; see Korn et al.
[\emph{J. Statist. Plann. Inference} \textbf{124} (2004) 379--398];
Romano and Wolf
(2007). However, the theoretical validity of the latter is still
largely unknown.
This paper provides a careful study of this heuristic: first, we extend
this approach by using a notion of
``bounding device'' that allows us to cover a wide range of critical
values, including those that adapt to $m_0$, the number of true null
hypotheses. Second, the theoretical validity of the latter is
investigated both nonasymptotically and asymptotically.
Third, we introduce suitable modifications of this heuristic that
provide new methods,
overcoming the existing procedures with a proven FDP control.
\end{abstract}



\end{frontmatter}

\section{Introduction}\label{sec1}

\subsection{Motivation} \label{sec:motiv}

Assessing significance in massive data is an important challenge of
contemporary statistics, which becomes especially difficult when the
underlying errors are
correlated. Pertaining to this class of high-\break dimensional problems, a
common issue is to make simultaneously a huge number $m$ of $0/1$
decisions with a valid control of the overall amount of false
discoveries (items declared to be wrongly significant). In this
context, a convenient way to account for false discoveries is the false
discovery proportion (FDP) that corresponds to the proportion of errors
among the items declared as significant (i.e., ``$1$'') by the procedure.

The Benjamini and Hochberg (BH) procedure has been widely popularized
after the celebrated paper \citet{BH1995} and is shown to control the
\textit{expectation} of the FDP, called the false discovery rate (FDR),
either theoretically under constrained dependency structures [see \citet
{BY2001}] or with simulations; see \citet{KW2008}. However, many authors
have noticed that the distribution of the FDP of BH procedure can be
affected by the dependencies [see, e.g., \citet
{Korn2004,DR2011,GHS2013}], which makes the use of the BH procedure
questionable.

\begin{figure}

\includegraphics{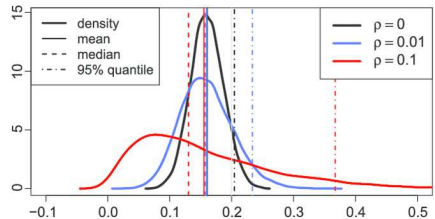}

\caption{Fitted density of the false discovery proportion of the BH
procedure when increasing the dependence.
$m=1000$, $m_0=800$ (number of true null hypotheses), $10^4$
simulations, Gaussian one-sided equicorrelated model.}
\label{FDPdep}
\end{figure}

To illustrate further this phenomenon, Figure~\ref{FDPdep} displays the
distribution of the FDP of the BH procedure in the classical one-sided
Gaussian multiple testing framework, when the $m$ test statistics are
all $\rho$-equicorrelated. As $\rho$ increases, the distribution of the
FDP becomes less concentrated and turns out to be drastically skewed
for $\rho=0.1$ (in particular, it falls outside the Gaussian regime).
Clearly, in this case, the mean fails to describe accurately the
overall behavior of the FDP distribution.
In particular, although the mean of the FDP is below $0.2$ [as proved
in \citet{BY2001}], the true value of FDP is not ensured to be small in
this case.

An alternative proposed in \citet{GW2004,PGVW2004,LR2005} is to control
the $(1-\zeta)$-quantile of the FDP distribution at level $\alpha$,
that is, to assert
%
%
\begin{equation}
\label{equFDPcontrolintro} \P(\FDP>\alpha)\leq\zeta.
\end{equation}
While taking $\zeta=1/2$ into \eqref{equFDPcontrolintro} provides a
control of the median of the FDP, taking $\zeta=0.05$ ensures that the
FDP does not exceed $\alpha$ with probability at least $95\%$.
Markedly, Figure~\ref{FDPdep} shows that the $(1-\zeta)$-quantiles of
the FDP distribution are substantially affected by the dependencies,
but not equally for all the $\zeta$'s: while the $95\%$-quantile
gets substantially larger, the median gets slightly smaller.
This suggests that the BH procedure is
much too optimistic for a $95\%$-quantile control, but is actually too
conservative for a FDP median control.
Overall, this reinforces the fact that in the presence of strong
dependence, controlling the $(1-\zeta)$-quantile of the FDP is an
essential task, not covered by the BH procedure.

%

\subsection{RW's heuristic and main contributions of this paper}

The problem of finding multiple testing procedures ensuring the control
\eqref{equFDPcontrolintro} has received growing attention in the last
decades; see, for instance,
\citet{CT2008},
\citet{DL2008},
\citet{GHS2013},
\citet{GR2007},
\citet{LR2005},
\citeauthor{RS2006b}
(\citeyear{RS2006,RS2006b}),
\citet{RW2007},
\citet{Roq2011},
\citet{RV2010}.
However, existing procedures with a proven FDP control are in general
too conservative. This increases the interest of simple and general
heuristics that work ``fairly.''
\citet{RW2007}, themselves referring to \citet{Korn2004}, have proposed
such a heuristic. It is called RW's heuristic in the sequel and
can be formulated as follows.

\textit{Start from a family $\mtc{R}_k$, $k\in\{1,\ldots,m\}$, of procedures
such that for all $k$, with probability at least $1-\zeta$, the
procedure $\mtc{R}_k$ makes less than $k-1$ false discoveries.
Then, choose some $\wh{k}$ such that $({\wh{k}}-1)/R_{{\hat{k}}}
\leq
\alpha$, where $R_k$ denotes the number of rejections of $\mtc{R}_k$.
Finally use $\mtc{R}_{\hat{k}}$.}

Note that, in the original formulation, $\wh{k}$ was constrained to be
chosen such that any $k'$ with $k'<\wh{k}$ should also satisfy
$({{k'}}-1)/R_{{{k'}}} \leq\alpha$ (``step-down'' approach). This
constraint is not necessarily applied here (e.g., ``step-up'' approach
is allowed).
The rationale behind this principle is that, for each $k$, the FDP of
$\mtc{R}_k$ is bounded by $(k-1)/R_k$ with probability $1-\zeta$, so
that the FDP of $\mtc{R}_{\hat{k}}$ should be smaller than $({\wh
{k}}-1)/R_{{\hat{k}}} \leq\alpha$ with probability $1-\zeta$, which
entails \eqref{equFDPcontrolintro}. However, as it is, this argument is
not rigorous because it does not take into account the fluctuations of
$\wh{k}$.

This heuristic has been theoretically justified (in the step-down form)
in settings where the $p$-values under the null are independent of the
$p$-values under the alternative [full independence in \citet{GR2007};
alternative $p$-values all equal to $0$ in \citet{RW2007}]. Since these
situations rely on an independence assumption, and since the FDP is
particularly interesting under dependence, it seems appropriate to study
the precise behavior of this method in ``simple'' dependent cases.
Thus our study is guided by the case where the dependencies are
\textit{known}, Gaussian multivariate or carried by latent variables.

In a nutshell, this paper makes the following main contributions:
\begin{itemize}[$-$]
\item[$-$] It provides a general framework in which RW's heuristic can be
investigated, by building the initial $k$-FWE critical values with
``bounding devices'': a strong interest is the possibility to build
critical values that ``adapt'' to $m_0$, the number of true nulls.
This allows to encompass many procedures, either new or previously known.
\item[$-$] We show that RW's heuristic may fail to control the FDP
nonasymptotically (even under its step-down form). Two corrections
that provably control the FDP are introduced. By using simulations, we
show that the resulting procedures are more powerful than those
previously existing.
\item[$-$]
We provide some asymptotic properties of RW's heuristic (in its step-up form):
first, we show that it is valid under weak dependence. In addition, we argue
that the interest of the latter is only moderate by proving that the simple
BH procedure is also valid in this case. Second, we provide particular types
of strong dependence for which RW's heuristic can be justified.
As a simple illustration, in a $\rho$-equicorrelated one-sided Gaussian
framework, we prove the asymptotic FDP control holds for
the step-up procedure using the following new critical values:
%
%
\begin{equation}
\label{simpletaul} \tau_{\ell}= \ol{\Phi} \bigl({\rho}^{1/2} {\ol{
\Phi}}^{-1} (\zeta) + (1-\rho )^{1/2} {\ol{
\Phi}}^{-1} (\alpha{\ell}/m) \bigr),\qquad 1\leq{\ell}\leq m,
\end{equation}
where $\ol{\Phi}$ is the upper-tail of the standard normal distribution.
\end{itemize}

\begin{figure}

\includegraphics{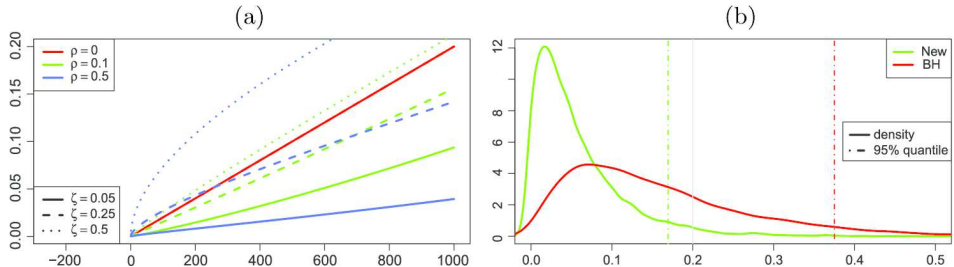}

\caption{\textup{(a)} plot of the critical values \protect\eqref{simpletaul}
in function of ${\ell}$. \textup{(b)} same as Figure~\protect\ref{FDPdep}
but only for $\rho=0.1$ and by adding the new step-up procedure using
\protect\eqref{simpletaul}. $m=1000$, $\alpha=0.2$,
Gaussian one-sided $\rho$-equi-correlated model.}
\label{FDPdepnew}
\end{figure}

Finally, let us emphasize that the critical values \eqref{simpletaul}
allow us to describe how the quantities $\alpha$, $\zeta$ and $\rho$
come into play when controlling (asymptotically) the FDP. Taking $\rho
=0$ just gives Simes's critical values, and thus the BH procedure,
whatever $\zeta$ is. The asymptotic FDP control can be explained in
this case by the fast concentration of the FDP of BH around its
expectation as $m$ grows to infinity under independence; see, for
example, \citet{Neu2008}.
Now, for $\rho>0$, the new critical values are markedly different from
the BH critical values: taking $\zeta=1/2$ leads to less conservative
critical values (if $\alpha\leq1/2$), while taking $\zeta$ smaller can
lead to more conservativeness (as expected); see Figure~\ref{FDPdepnew}(a)
for an illustration. Finally, we plot in Figure~\ref{FDPdepnew}(b)
the density of the FDP of the step-up procedure using the new critical
values \eqref{simpletaul} for $\zeta=0.05$. As one might expect,
compared to the BH procedure, the density has been shifted to the left
so that the $95\%$-quantile of the FDP of the novel procedure is below
$\alpha$.

%

\subsection{Multiple testing framework}\label{sec:setting}

We observe a random variable $X$, whose distribution belongs to some
set $\mtc{P}$. For $m\geq2$,
we define a setting for performing $m$ tests simultaneously by
introducing a true/false null parameter $H\in\{\mbox{0,1}\}^m$ and a set of
associated distributions $\mtc{P}_H\subset\mtc{P}$
which are candidates to be the distribution of $X$ under the
configuration $H$. We denote ${\mtc{H}}_0(H)=\{i\dvtx H_i=0\}$,
$m_0(H)=\sum_{i=1}^m (1-H_i)$ and ${\mtc{H}}_1(H)=\{i\dvtx H_i=1\}$, $m_1(H)=\sum_{i=1}^m
H_i$ the set/number of true and false nulls, respectively. The basic
assumption is the following: for all $i\in\{1,\ldots,m\}$, there is a
$p$-value $p_i(X)$ satisfying the following assumption:
%
\[
\forall H \in\{0,1\}^m \mbox{ with } H_i=0, \forall P\in\mtc{P}_H, \forall t\in[0,1],\qquad \P_{X\sim P}(p_i(X)\leq t)\leq t.\vadjust{\goodbreak}
\]
In this paper, a leading example is the one-sided location model,
%
%
\begin{equation}
\label{equonesided}X_i=H_i \mu_i+
Y_i, \qquad 1\leq i \leq m,
\end{equation}
where $H\in\{0,1\}^m$, $\mu\in(\R_+\setminus\{0\})^m$ and $Y$ is a
$m$-dimensional centered random vector with identically distributed components.
Then the $p$-values are given by $p_i(X)=\ol{F}(X_i)$, where $\ol{F}
(x)=\P(Y_1\geq x)$, $x\in\R$.
Note that this model implicitly assumes that the $p$-values under the
null are uniformly distributed.
In this paper, we will often assume that the joint distribution of the
noise $Y$ is known, and we consider the two following models for $Y$:
\begin{itemize}[$-$]
\item[$-$] Gaussian: $Y$ is a Gaussian vector with covariance matrix
$\Gamma$ (such that $\Gamma_{i,i}=1$ for simplicity),
in which case $Y_1\sim\mathcal{N}(0,1)$ and  $\Fbar$ is denoted by $\pnorm$.
A simple particular  case is the equi-correlated case,
%
%
\renewcommand{\theequation}{Gauss-$\rho$-equi}
\begin{equation}\qquad
\label{equi-correlated}
\Gamma_{i,j}=\rho\qquad \mbox{for all } i\neq j,
\mbox{ where } \rho \in \bigl[-(m-1)^{-1},1\bigr].
\end{equation}

\item[$-$] Mixture of ($1$-)factor models: the distribution of $Y$ is
given by
%
%
\renewcommand{\theequation}{facmod}
\begin{equation}
\label{facmodel}
Y_i=c_i W +
\xi_i, \qquad 1\leq i \leq m,
\end{equation}
where $c_i$, $1\leq i \leq m$, are i.i.d., $\xi_i$, $1\leq i \leq m$,
are i.i.d., $W$ is a random variable and $(c_i)_{1\leq i\leq m}$, $(\xi
_i)_{1\leq i\leq m}$ and $W$ are independent.
Also, the distributions of $W$, $c_1$ and $\xi_1$ are assumed to be
known, so that the function $\ol{F}$ is known, and the $p$-values can
be computed. A simple particular case is obtained as follows:
for $\rho\in[0,1]$,
%
%
\renewcommand{\theequation}{alt-$\rho$-equi}
\begin{equation}
\label{alt-equi-correlated}
Y_i= \varepsilon_i
\rho^{1/2} W + (1-\rho)^{1/2} \zeta_i,\qquad 1\leq i
\leq m,
\end{equation}
where $W,\zeta_1, \ldots, \zeta_m$ are i.i.d. $\mathcal{N}(0,1)$
and are
independent of $\varepsilon_1,\ldots,\varepsilon_m$ which are i.i.d.
random signs following the distribution $(1-a) \delta_{-1}+a \delta_1$,
for a parameter $a\in[0,1]$.
\end{itemize}
While the Gaussian model is classical and widely used, \eqref{facmodel}
is useful to model a strong dependence, through the factor $W$.
When the $c_i$'s are deterministic, the latter is often referred to as
a \textit{one factor model} in the literature, see, for example, \citet
{LS2008,FKC2009,FHG2012}. Here, the $c_i$'s are unknown and taken
randomly with a prescribed distribution.
From an intuitive point of view, \eqref{facmodel} is modeling
situations where some of the measurements have been deteriorated by
unknown nuisance factors $c_i W$, $1\leq i\leq m$. For instance,
choosing $c_i\in\{0,1\}$ corresponds to simultaneously deteriorate the
measurements of some unknown sub-group $\{1\leq i\leq m\dvtx c_i=1\}
\subset
\{1,\ldots,m\}$.
Furthermore, note that while the model \eqref{alt-equi-correlated}
covers \eqref{equi-correlated} when $\rho\geq0$ by taking $a=0$,
\eqref{alt-equi-correlated} is able to include negative dependence
between some of the $Y_i$'s.

In \eqref{facmodel}, a quantity of interest throughout the paper is the
probability that a $p$-value is below $t$ conditionally on $W=w$ (under
the null).\vadjust{\goodbreak} According to the particular setting that is at hand, this
probability can be written as follows: for $\rho\in[0,1)$, $w\in\R$,
%
%
\renewcommand{\theequation}{$F_0$-facmod}
\begin{equation}\qquad
\label{Fzero-facmodel}
\mbox{\hspace*{64pt}$F_{0}(t,w)=\P\bigl(\ol{F}(c_1w+\xi_1)\leq
t\bigr)=\E \bigl[ \ol{F}_\xi \bigl(\ol {F}^{-1}(t)-c_1w
\bigr) \bigr];$}\vspace*{-12pt}
\end{equation}
\renewcommand{\theequation}{$F_0$-alt-$\rho$-equi}
\begin{equation}
\label{Fzero-alt-equi-correlated}
\mbox{\hspace*{34pt}$F_{0}(t,w)= (1-a) f(t,-w,\rho) +a f(t,w,\rho);$}\vspace*{-12pt}
\end{equation}
\renewcommand{\theequation}{$F_0$-Gauss-$\rho$-equi}
\begin{equation}
\label{Fzero-equi-correlated}
\mbox{\hspace*{-70pt}$F_{0}(t,w)= f(t,w,\rho),$}
\end{equation}
where $\ol{F}_\xi(x)=\P(\xi_1\geq x)$ and
$
f(t,w,\rho)= \ol{\Phi} (({{\ol{\Phi}}^{-1} (t) - {\rho
}^{1/2} w)}/{(1-\rho
)^{1/2}} )$.

%
\begin{remark}[(Modifications of the test statistics)]
Let us consider the model \eqref{facmodel}, where $c_i$ is equal to
some known constant; \eqref{equi-correlated} is one typical instance.
Then, as noted by a referee, applying a re-centering operation to the
$X_i$'s makes the
factor $W$ disappear,
and thus can lead to better test statistics
(if the bias induced  by this operation is not too large);
see Section~S-1 in the \hyperref[suppA]{Supplementary Material} for more details on this issue.
In this respect, our work is particularly relevant in cases where $W$
cannot be estimated (but has a known distribution).
On the other hand, we believe that model \eqref{equi-correlated} keeps
the virtue of simplicity and hence remains interesting when studying
procedures that are supposed to deal with strong dependencies.
Hence while our procedures will in general not be restricted to model
\eqref{equi-correlated}, we will also use this model for illustrative
purposes throughout the paper.
\end{remark}

In the Gaussian case, the joint distribution of the $p$-values under
the null $(p_i, i\in{\mtc{H}}_0(H))$ depends, in general, on the subset
${\mtc{H}}
_0(H)$. Obviously, in this case, we do not want to explore the $m
\choose m_0(H)$ possible subsets of $\{1,\ldots,m\}$ in our inference,
which inevitably should arise when our procedure fits to such a
dependence structure. To circumvent this technical difficulty, we can
add random effects to our model. This makes $H$ become random.
More formally, we distinguish between the two following models:
\begin{itemize}
\item[$-$] Fixed mixture model: the parameter $H$ is fixed by advance and
unknown. Overall, the parameters of the model are given by $\theta
=(H,P)$ to be chosen in the set
\[
\Theta^F=\bigl\{(H,P)\dvtx H\in\{0,1\}^m,P\in
\mtc{P}_H\bigr\}.
\]
\item[$-$] Uniform mixture model: the number of true null $m_0\in\{
0,1,\ldots,m\}$ is unknown and fixed by advance, while $H$ is a random
vector distributed in such a way that ${\mtc{H}}_0(H)$ is randomly generated
(independently and\break previously of the other variables), uniformly in the
subsets of $\{1,\ldots,m\}$\break of cardinal $m_0$. The parameters of the
model are given by $\theta=\break (m_0,(P_H)_{H\dvtx m_0(H)=m_0})$, to be
chosen in
the set
\begin{eqnarray*}
&&\Theta^{U}=\bigl\{(m_0,(P_H)_{H:m_0(H)=m_0} )\telque m_0\in \{0,1,\ldots,m\} ,\\
&&\hspace*{80pt} P_H \in \mtc{P}_H \mbox{for all $H:m_0(H)=m_0$}\bigr\}.
\end{eqnarray*}
In this model, the distribution of $X$ conditionally on $H$ is $P_H$.
\end{itemize}
While the fixed mixture model is the most commonly used model for
multiple testing, the
uniform mixture model
is new to our knowledge and follows the general philosophy of models
with random effects; see \citet{ETST2001}. It is convenient for the
adaptation issue w.r.t. $m_0$, as we will see later on.
With some abuse, we denote $m_0(\theta)$, ${\mtc{H}}_0(\theta)$ (or $m_0$,
${\mtc{H}}_0$ when not ambiguous) the number of true nulls in the
fixed/uniform mixture models.
In the sequel, $\Theta$ denotes either $\Theta^F$ or $\Theta^{U}$.

\subsection{Type I error rates}
First, for $t\in[0,1]$, denote by $V_m(t)=\sum_{i=1}^m (1-H_i){\mbf
{1}\{p_i(X)\leq t\}}$ and $R_m(t)=\sum_{i=1}^m{\mbf{1}\{p_i(X)\leq t\}
}$ the number
of false discoveries and the number of discoveries (at threshold $t$),
respectively.
For some pre-specified $k\in\{1,\ldots,m\}$ and some thresholding method
$\wh{t}_{m}\in[0,1]$ (potentially depending on the data), the
$k$-family-wise error rate ($k$-FWER) is defined as the probability
that more than $k$ true nulls have a $p$-value smaller than $\wh
{t}_{m}$; see, for example, \citet{Hom1988,LR2005}. Formally, for
$\theta
\in\Theta$ (in one of the models defined in Section~\ref{sec:setting}
and $\Theta$ being the corresponding parameter space),
\renewcommand{\theequation}{\arabic{equation}}
\setcounter{equation}{3}
\begin{equation}
\label{kFWERformula} \mbox{$k$-FWER}(\wh{t}_{m})=\P_\theta
\bigl(V_m(\wh{t}_{m})\geq k\bigr).
\end{equation}
Note that $k=1$ corresponds to the traditional family-wise error rate
(FWER). From \eqref{kFWERformula}, providing $\mbox{$k$-FWER}(\wh
{t}_{m})\leq\zeta$
(for all $\theta\in \Theta$), ensures that, with probability at least $1-\zeta$, less than
$k-1$ false discoveries are made by the thresholding procedure $\wh{t}_{m}$.

Next, for some threshold $t\in[0,1]$, define the false discovery
proportion at threshold $t$ as follows:
%
%
\begin{equation}
\label{FDP} \FDP_m(t)= \frac{V_m(t)}{R_m(t)\vee1}.
\end{equation}
Note that the quantity $\FDP_m(t)$ is random and not observable because
it depends on the unknown process $V_m(t)$. Controlling the FDP via a
threshold $t=\wh{t}_{m}$ (potentially depending on the data)
corresponds to the following probabilistic bound:
%
%
\begin{equation}
\label{FDPcontrol} \forall\theta\in\Theta\qquad \P_{\theta} \bigl(\FDP_m
(\wh {t}_{m} )\leq\alpha \bigr)\geq1-\zeta,
\end{equation}
for some pre-specified values $\alpha,\zeta\in(0,1)$. As mentioned
before, \eqref{FDPcontrol} corresponds to upper-bounding the $(1-\zeta
)$-quantile of the distribution of\break $\FDP_m (\wh{t}_{m} )$ by
$\alpha$.
Since $\FDP_m(t)>\alpha$ is equivalent to $V_m(t)\geq\lfloor\alpha
R_m(t)\rfloor+1$, the FDP control and the $k$-FWER control are
intrinsically linked.

From a historical point of view, the introduction of the FDP goes back
to Eklund in the 1960s [as reported in \citet{See1968}], who has
presented the FDP as a solution to the ``mass-significance problem.''
Much later, the seminal paper of \citet{BH1995} has widely popularized
the use of the FDP in practical problems by introducing and studying
the false discovery rate (FDR), which corresponds to the expectation of
the FDP.

\subsection{Step-up and step-down procedures}

Let us consider the ordered $p$-values $p_{(1)}\leq\cdots\leq p_{(m)}$.
Consider a nondecreasing sequence $(\tau_{\ell})_{1\leq{\ell}\leq
m}$ of
nonnegative values, referred to as the \textit{critical values}. The
corresponding step-up (resp., step-down) procedure is defined as
rejecting the $p$-values smaller than $\tau_{\hat{{\ell}}}$, where
$\wh
{{\ell}
}$ is defined by either of the two following quantities (with the
convention $p_{(0)}=0,\tau_0=0$):
%
%
\renewcommand{\theequation}{SU}
\begin{equation}
\max\bigl\{{\ell}\in\{0,1,\ldots,m\}\mbox{ such that } p_{({\ell})} \leq
\tau _{\ell}\bigr\} \label{SU-classic};\vspace*{-12pt}
\end{equation}
\renewcommand{\theequation}{SD}
\begin{equation}
\max\bigl\{{\ell}\in\{0,1,\ldots,m\}\mbox{ such that } \forall{\ell
}'\in\{ 0,1,\ldots ,{\ell}\}, p_{({\ell}')} \leq
\tau_{{\ell}'} \bigr\}.\label{SD-classic}%
\end{equation}
Let us also recall the so-called \textit{switching relation}:
$p_{({\ell})}
\leq\tau_{\ell}$ is equivalent to $R_m(\tau_{{\ell}}) \geq{\ell}$.
This entails $R_m(\tau_{\hat{\ell}})=\wh{{\ell}}$ both in the
step-up and
step-down cases.

\section{Building $k$-FWE-based critical values}\label{seckFWERcrit}

\subsection{Revisiting RW's heuristic}
Starting from arbitrary critical values $(\tau_{\ell})_{1\leq{\ell
}\leq m}$,
and by taking an integer $\wh{{\ell}}$ such that $R_m(\tau_{\hat
{{\ell}}}) =
{\wh{{\ell}}}$, we have
%
\renewcommand{\theequation}{\arabic{equation}}
\setcounter{equation}{6}
\begin{eqnarray}
\label{heuristic} \P_\theta\bigl( \FDP_m(\tau_{\hat{\ell}})>
\alpha\bigr) &=&\P_\theta \bigl(V_m(\tau _{\hat{\ell}
}) >
\alpha R_m(\tau_{\hat{\ell}}) \bigr)
\nonumber
\\[-8pt]
\\[-8pt]
\nonumber
&=&\P_\theta\bigl( V_m(\tau_{\hat{\ell}}) \geq\lfloor
\alpha\wh{{\ell} }\rfloor+1\bigr).
\end{eqnarray}
Hence, by taking $\tau_{\ell}$ such that $( \lfloor\alpha{{\ell
}}\rfloor+1
)\mbox{-FWER}(\tau_{\ell}) \leq\zeta$ for all ${\ell}$, we should
get that
\eqref{heuristic} is below $\zeta$. However, as already mentioned, the
above reasoning does not rigorously establish \eqref{FDPcontrol} (with
$\wh{t}_{m}=\tau_{\hat{\ell}}$) because it implicitly assumes that
$\wh
{{\ell}
}$ is deterministic.
Nevertheless, this heuristic is a suitable starting point for building
critical values related to the FDP control.

\subsection{Bounding device}\label{sec:nullquantile}

Let us consider either the fixed model $\Theta=\Theta^F$ or the uniform
model $\Theta=\Theta^U$.
First, let us define a bounding device as any function $B^0_m \dvtx
(t,k,u)\mapsto B^0_m(t,k,u)\in[0,1]$, defined for $t\in[0,1]$, $k\in
\{
1,\ldots,m\}$ and $u\in\{0,\ldots,m\}$, which is nonincreasing in $k$,
with $B^0_m(0,k,u)=0$ for all $u,k$, $B^0_m(t,k,u)=0$ for all $t\in
[0,1]$ whenever $u<k$, and such that for all $t\in[0,1]$, $k\in\{
1,\ldots,m\}$ and $u\in\{k,\ldots,m\}$, we have
%
%
\renewcommand{\theequation}{Bound}
\begin{equation}
B^0_m(t,k,u)\geq\mathop{\sup_{\theta\in\Theta}}_{ m_0(\theta
)=u}
\bigl\{\P_{\theta}\bigl( V_m(t) \geq k\bigr) \bigr\}
\label{boundingfunction}.%
\end{equation}
Now, define for $t\in[0,1]$, $k\in\{1,\ldots,m\}$ and ${\ell}\in\{
k,\ldots
,m\}
$, the quantities
%
%
\renewcommand{\theequation}{Bound-nonadapt}
\begin{equation}
\ol{B}_m(t,k) = \sup_{0\leq u \leq m} \bigl\{
B^0_m(t,k,u) \bigr\}; \label{boundingfunctionnonadapt}\vspace*{-12pt}
\end{equation}
\renewcommand{\theequation}{Bound-adapt}
\begin{equation}
\hspace*{48pt}\wt{B}_m(t,k,{\ell}) = \sup_{k\leq k'\leq{\ell}}
\Bigl\{\sup_{
0\leq u
\leq
m-{\ell}+k'} B^0_m
\bigl(t,k',u\bigr) \Bigr\},\label{boundingfunctionadapt} %
\end{equation}
which are additionally assumed to be nondecreasing and left-continuous
in $t$.
Note that $\ol{B}_m(t,k)$ and $\wt{B}_m(t,k,{\ell})$ are both
nonincreasing in $k$.

%
\begin{definition}\label{defkFWE}
Let us consider a bounding device $B^0_m(t,k,u)$ and the above
associated quantities $\ol{B}_m(t,k)$ and $\wt{B}_m(t,k,{\ell})$.
Then the nonadaptive (resp., adaptive, oracle) $k$-FWE-based critical
values associated to the bounding function $B^0_m$ are defined as
follows (resp.):
%
%
\renewcommand{\theequation}{\arabic{equation}}
\setcounter{equation}{7}
\begin{eqnarray}
\ol{\tau}_{\ell}&=& \max \bigl\{ t\in[0,1] \dvtx \ol{B}_m\bigl(t,
\lfloor \alpha {\ell}\rfloor+1\bigr) \leq\zeta \bigr\},\qquad 1\leq{\ell}\leq
m;\label
{equtauknonadapt}
\\
\wt{\tau}_{\ell}&= &\max \bigl\{ t\in[0,1] \dvtx \wt{B}_m\bigl(t,
\lfloor \alpha {\ell}\rfloor+1,{\ell}\bigr) \leq\zeta \bigr\},\qquad 1\leq{\ell}\leq
m;\label
{equtaukadapt}
\\
{\tau}^0_{\ell}&=& \max \bigl\{ t\in[0,1] \dvtx
{B}^0_m\bigl(t, \lfloor \alpha {\ell}
\rfloor+1,m_0\bigr) \leq\zeta \bigr\},\qquad 1\leq{\ell}\leq
m.\label
{equtaukoracle}
\end{eqnarray}
\end{definition}

The above definition implies that $(\ol{\tau}_{\ell})_{1\leq{\ell
}\leq m}$,
$(\wt{\tau}_{\ell})_{1\leq{\ell}\leq m}$ and $({\tau}^0_{\ell
})_{1\leq{\ell}
\leq m}$
are nondecreasing sequences, so that they can be used as critical
values. The critical values $\wt{\tau}_{\ell}$, ${\ell}=1,\ldots
,m$, are
said to
be \textit{adaptive} because they implicitly (over-)estimate $m_0$ by
%
%
\renewcommand{\theequation}{\arabic{equation}}
\setcounter{equation}{10}
\begin{equation}
\label{equml} m({\ell})=m-{\ell}+\lfloor\alpha{\ell}\rfloor+1.
\end{equation}
In the literature, this way to adapt to $\pi_0$ is often referred to as
\textit{one-stage} [in contrast to \textit{two-stage}; see \citet
{BKY2006,Sar2008,BR2009}]. It has been proved to be asymptotically
optimal in a specific sense; see \citet{FDR2009}.
Also, $\ol{\tau}_{\ell}\leq\wt{\tau}_{\ell}$ for all ${\ell}$;
that is,
adaptation always leads to less conservative critical values.
Finally, it is worth to check that
$\ol{\tau}_m\leq\wt{\tau}_m <1$ [this comes from $B_m^0(1,k,u)=1$ for
all $ u\geq k$] so that the output $\wh{{\ell}}$ of the step-up algorithm
is not identically equal to $m$.

%
\begin{figure}

\includegraphics{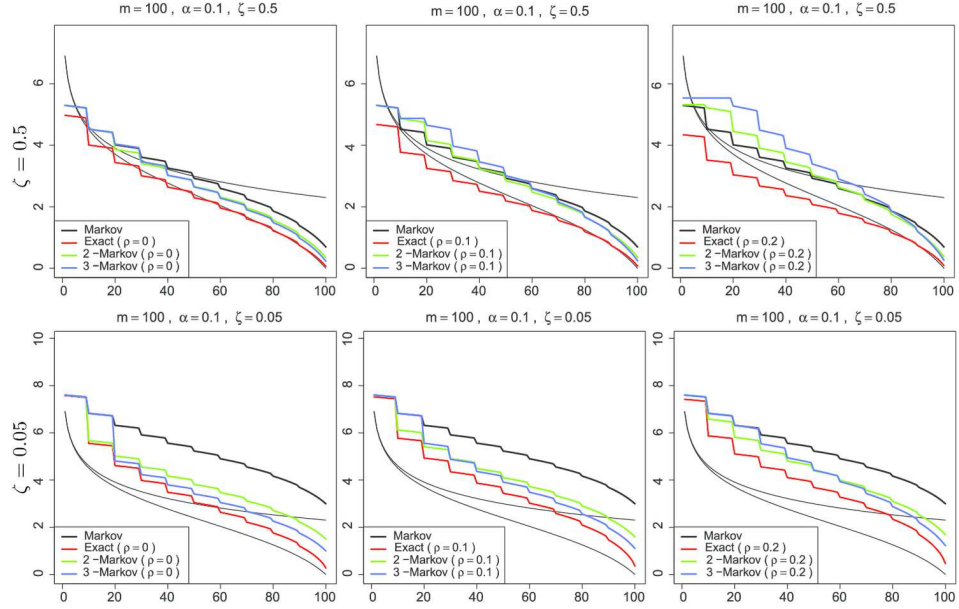}

\caption{Plot of $-\log(\wt{\tau}_{\ell})$ in function of ${\ell
}\in\{
1,\ldots
,m\}$, for $k$-FWE-based critical
values obtained with several types of bounding devices and assuming
\protect\eqref{equi-correlated}. For comparison,
the solid thin black line corresponds either to the BH critical values
$-\log(\alpha{\ell}/m)$ (nonadaptive)
or the AORC critical values $-\log(\alpha{\ell}/(m-(1-\alpha){\ell
}))$ defined
in \citet{FDR2009} (adaptive).}
\label{critvalues2}
\end{figure}

\subsection{Examples}\label{sec:examples}

We provide below three examples of bounding devices: Markov, $K$-Markov
and Exact. Instances of resulting critical values are displayed in
Figure~\ref{critvalues2} under Gaussian equi-correlation (see Figure
S-2 
for similar pictures under alternate equi-correlation). As we will see,
while the exact bounding device leads to the largest critical values,
the Markov-type devices are still useful because they can offer finite
sample controls.
Also note that in all these examples, we have $\wt{B}_m(t,k,{\ell}
)=B_m^0(t,k,m-{\ell}+k)$.

\textit{Markov.}
By Markov's inequality, we have
%
%
\begin{equation}
\P_\theta\bigl(V_m(t)\geq k\bigr) \leq\frac{\E_\theta(V_m(t))}{k} =
\frac{m_0
t}{k}=:B^0_m(t,k,m_0).\label{MarkovDevice}
\end{equation}
Since $\ol{B}_m(t,k)=mt/k$ and $\wt{B}_m(t,k,{\ell}) = (m-{\ell
}+k)t/k$, this
gives rise to the critical values
%
%
\begin{equation}
\ol{\tau}_{\ell}= \frac{\zeta(\lfloor\alpha{\ell}\rfloor
+1)}{m} ;\qquad \wt {\tau }_{\ell}=
\frac{\zeta(\lfloor\alpha{\ell}\rfloor+1)}{m({\ell
})},\qquad 1\leq{\ell} \leq m,\label{criticalvaluesMarkov}
\end{equation}
where $m({\ell})$ is defined by \eqref{equml}.
The adaptive critical values $(\wt{\tau}_{\ell})_{1\leq{\ell}\leq
m}$ are
those proposed by \citet{LR2005}.
Note that these critical values do not adapt to the underlying
dependence structure of the $p$-values.

%
%
%
%
%
%
%
%
%
%
%

\textit{$K$-Markov.}
When $\zeta$ is small, the Markov device can be too conservative, and
we might want to use a sharper tool.
Let $K\geq1$ be an integer. As suggested in \citet{GHS2013} (for
$K=2$), we can use the following bound: for $k\geq K$,
%
%
\begin{equation}
\label{equ:majKmarkov} {\mbf{1}\bigl\{V_m(t) \geq k\bigr\}}\leq
\frac{1}{ {k \choose K} }\sum_{X\subset{\mtc{H}}
_0\dvtx|X|=K} {\mbf{1}\Bigl\{\max
_{ i\in X}\{ p_i\}\leq t\Bigr\}},
\end{equation}
which leads to upper-bounding $ \sup\{ \P_\theta(V_m(t)\geq
k),\theta
\in\Theta, m_0(\theta)=u\}$ by
%
%
\begin{equation}
B_m^0(t,k,u)=\frac{u(u-1)\cdots(u-K+1)}{k(k-1)\cdots(k-K+1)} \mathop {\sup
_{\theta\in\Theta}}_{ m_0(\theta)=u} \Bigl\{ P_\theta \Bigl(\max
_{
i\in X_0}\{ p_i\}\leq t \Bigr) \Bigr\} ,\label{KMarkovbounding}
\end{equation}
where $X_0\subset{\mtc{H}}_0$ with $|X_0|=K$ is let arbitrary. In the
latter, we implicitly assume that for all $\theta$, the probability
$P_\theta(\max_{ i\in X}\{ p_i\}\leq t)$ depends on $X$ only through $|X|=K$.
When $k<K$, bound \eqref{KMarkovbounding} is useless in essence, so we
replace it by the simple Markov device, by letting $B^0_m(t,k,u)
=(ut/k) \vee B_m^0(t,K,m)$ if $1\leq k \leq K-1$. Note that the
operator ``$\vee$'' in the last display is added to keep the
nonincreasing property w.r.t. $k$.
As a first illustration, in the one-sided location model \eqref
{equonesided}, if
%
%
\renewcommand{\theequation}{\mbox{Indep}}
\begin{equation}\qquad
\label{indep}
\mbox{$p_i,i\in{\mtc{H}}_0$,
are mutually independent } \bigl(\mbox{cond. on } H\mbox{ in model }\Theta^U\bigr),
\end{equation}
we have $P_\theta(\max_{ i\in X_0}\{ p_i\}\leq t)=t^K$,
which entails
\[
(\wt{\tau}_{\ell})^K=\zeta\frac{(\lfloor\alpha
{\ell}
\rfloor
+1)(\lfloor\alpha{\ell}\rfloor+1-1)\cdots(\lfloor\alpha{\ell
}\rfloor
+1-K+1)}{m({\ell})(m({\ell})-1)\cdots(m({\ell})-K+1)}
\]
for
${\ell}\geq\lceil(K-1)/\alpha\rceil$.
A second illustration is the one-factor model \eqref{facmodel}, for
which $P_\theta(\max_{ i\in X_0}\{ p_i\}\leq t)= \E
[F_{0}(t,W)^{K} ] $ where $F_{0}$ is defined by\break \eqref{Fzero-facmodel}.
Hence, inverting $B^0_m(t,\lfloor\alpha{\ell}\rfloor+1,m({\ell
}))=\zeta$
gives rise to critical values $\ol{\tau}_{\ell}$ and $\wt{\tau
}_{\ell}$, which
both take into account the dependence induced by the common factor $W$.

\textit{Exact.}
In some cases, closed-formulas can be derived for the RHS of \eqref
{boundingfunction}.
First, by assuming \eqref{indep},
the distribution of $V_m(t)$ is a binomial with parameters $(u,t)$. Hence,
$B_m^0(t,k,u)= \sum_{j=k}^{u} {u \choose j} t^{j} (1-t)^{u-j}$.
The corresponding adaptive critical values can be obtained by a
numerical inversion [these critical values were already proposed in
\citet{GR2007}]. Second, the following exact formula can be used in
model \eqref{facmodel}:
%
%
\renewcommand{\theequation}{\arabic{equation}}
\setcounter{equation}{15}
\begin{equation}\label{equ-multirhoExact}
{B}^0_m(t,k,u)= \E_W \Biggl[ \sum
_{j=k}^{u} \pmatrix{u \cr j} \bigl(F_{0}(t,W)
\bigr)^{j} \bigl(1-F_{0}(t,W)\bigr)^{u-j} \Biggr],
\end{equation}
where $F_0$ is defined by \eqref{Fzero-facmodel}.
Third, in a more general manner, nonadaptive threshold can be obtained
in the one-sided location model \eqref{equonesided}, provided that the
full joint distribution of $Y$ is known: for this, observe that
$V_m(t)$ is upper-bounded by the full null process
%
%
\begin{equation}
\label{fullnullprocess} V_m'(t)=m^{-1}\sum
_{i=1}^m {\mbf{1}\bigl\{\ol{\Phi}(Y_i)
\leq t\bigr\}},
\end{equation}
whose distribution can be approximated by a Monte Carlo method.
Finally, to obtain an adaptive threshold, we can make use of the
uniform model $\Theta^U$: the added random effect on ${\mtc{H}}_0$ entails
that $V_m(t)$ can be easily generated for each value of $u=m_0(\theta
)$. This leads to nonadaptive and adaptive critical values that
incorporate any pre-specified joint distribution for the noise $Y$
(e.g., $\Gamma$ in the Gaussian case).

\section{Finite sample results}\label{sec:finitebound}

\subsection{Preliminary results}

The following theorem gathers the only existing cases where RW's
heuristic has been proved to provide FDP control (to our knowledge).

%
\begin{proposition}[{[\citet{RW2007,GR2007,GHS2013}]}]\label{RWalreadyproved}
Consider some bounding device $B_m^0$ and the associated $k$-FWE-based
critical values $(\tau_{\ell})_{1\leq{\ell}\leq m}$, being either
adaptive or
not and computed either in the fixed mixture model ($\Theta=\Theta^F$)
or in the uniform mixture model ($\Theta=\Theta^U$). Let us consider
the corresponding number of rejections $\wh{{\ell}}$ of the associated
step-down \eqref{SD-classic} or step-up \eqref{SU-classic} procedure.
Then the FDP control \eqref{FDPcontrol} holds (with $\wh{t}_{m}=\tau
_{\hat{\ell}}$) in the following cases:
\begin{longlist}[(ii)]
\item[(i)] step-down algorithm and the null $p$-values $(p_i, i\dvtx H_i=0)$
are independent of the alternative $p$-values $(p_i, i\dvtx H_i=1)$;
\item[(ii)] step-down or
step-up algorithm with the Lehmann--Romano critical values, that is,
with $(\tau_{\ell})_{1\leq{\ell}\leq m}$ given by \eqref
{criticalvaluesMarkov}, and assuming that Simes's inequality is valid,
%
%
\begin{equation}
\label{equsimes}
\forall \theta\in\Theta\qquad \P_\theta\Biggl( \bigcup_{k= 1}^{m_0} \{q_{({k})} \leq \zeta k/m_0\} \Biggr)\leq \zeta,
\end{equation}
where $q_{(1)}\leq q_{(2)}\leq\cdots\leq q_{(m_0)}$ denote the ordered
$p$-values under the null.
\end{longlist}
\end{proposition}

Case (i) comes from inequalities established in \citet{LR2005,RW2007},
that we recall in Section~\ref{sec:unif} under an unified form; see
also Theorem~5.2 in \citet{Roq2011}. Note that it contains the case
where all the $p$-values under the alternative are equal to zero (Dirac
configuration). Case (ii) has been solved more recently in \citet
{GHS2013}. Here, it can be seen as a consequence of the following
general inequality; see Section~\ref{sec:prooffinitebound} for a proof.

%
\begin{proposition}\label{prop:finitebound}
Consider the setting of Proposition~\ref{RWalreadyproved} in the
step-down or step-up case.
Then we have for all $\theta\in\Theta$
%
%
\begin{eqnarray}
\label{equ-finiteboundSU} \P_\theta \bigl(\FDP_m(\tau_{\hat{{\ell}}}) >
\alpha \bigr) &\leq&\P_\theta \bigl( V_m\bigl(
\nu^0_{\hat{k}}\bigr) \geq\wh{k}\geq 1 \bigr)
\nonumber
\\[-8pt]
\\[-8pt]
\nonumber
&=&\P _\theta
\bigl( q_{(\hat{k})} \leq\nu^0_{\hat{k}}, \wh{k}\geq1 \bigr),
\end{eqnarray}
for $\wh{k}=V_m(\tau_{\hat{{\ell}}})$,
where
%
%
\begin{equation}
\nu^0_k=\max\bigl\{t\in[0,1]\dvtx {B}^0_{m}(t,{k},m_0)
\leq\zeta\bigr\} .\label
{oracle-critvalues}
\end{equation}
\end{proposition}

Proposition~\ref{RWalreadyproved}(ii) thus follows from
Proposition~\ref{prop:finitebound}, used with the adaptive Markov
bounding device; see \eqref{MarkovDevice}.
Markedly, Proposition~\ref{prop:finitebound} establishes that the FDP
control for adaptive $k$-FWE based critical values is linked to a
specific inequality between the null $p$-values and the bounding device
using the true value of $m_0$.

Further note that \eqref{equ-finiteboundSU} in Proposition~\ref
{prop:finitebound} is sharp whenever $m_0(\theta)=m$: in this case, the
LHS and RHS are both equal to the probability that $\wh{k}$($=\wh
{{\ell}}$)
is not zero, that is, that at least one ${\ell}\in\{1,\ldots,m\}$ is such
that $p_{({\ell})}\leq{\tau}_{{\ell}}$. For instance, in the
independent case
and using the exact device \eqref{equ-multirhoExact}, when $m=m_0=2$
and $\alpha=0.5$, we have $\ol{\tau}_{1}=\wt{\tau}_{1}=1-(1-\zeta
)^{1/2}$ and $\ol{\tau}_{2}=\wt{\tau}_{2}=\zeta^{1/2}$ and
$\P_\theta (\FDP_m(\tau_{\hat{{\ell}}}) >\alpha )=
2\zeta
-(1-(1-\zeta)^{1/2})(2\zeta^{1/2}-1 +(1-\zeta)^{1/2} )$.
We merely check that the latter is larger than $\zeta$ for all $\zeta
\in
(0,1)$. Also, simulations of Section S-5 
in the \hyperref[suppA]{Supplementary Material} [\citet{DR2014supp}] indicate that this exceeding can hold for a larger value
of $m$.
This establishes the following:\vadjust{\goodbreak}
%

\begin{fact}\label{factoracleSU}
RW's heuristic does not always provide a valid FDP control for finite
$m$ in its step-up form under independence.
\end{fact}

Now, an important question is to know whether RW's heuristic always
provides a valid FDP control for finite $m$ in its \textit{step-down} form.
First, we can merely check that the following cases can be added in
Proposition~\ref{RWalreadyproved} in the step-down case:
\begin{longlist}[(iii)]
\item[(iii)] for all $\theta\in\Theta$, $\lfloor\alpha b_\alpha
(m_0(\theta))\rfloor= 0$ (e.g., $m_0(\theta)\in\{1,m\}$ or $\lfloor
\alpha{m} \rfloor=0$);
\item[(iv)] under \eqref{equi-correlated} when $\rho=1$.
\end{longlist}
Note that (iii) contains the case $m_0=m$ which is problematic in the
step-up case.
A consequence is that any configuration for which the FDP control fails
should be searched outside cases (i), (ii), (iii) and (iv).
As a matter of fact, we found a numerical example under
equi-correlation when using the critical values $({\tau}^0_{\ell
})_{\ell}$
defined by \eqref{equtaukoracle}, with the exact device. To this end,
we have evaluated the exceedance probability of the FDP by the
\textit{exact calculations} proposed in \citet{RV2010,BDRV2014}. This
method is
time consuming for large $m$ but avoids the undesirable fluctuations
due to the Monte Carlo approximation while performing simulations.
Precisely, in model \eqref{equi-correlated}, when $m=30$, $\alpha=0.2$,
$\zeta=0.05$, $\rho=0.3$, $m_0=15$, $\mu_i=1.5$, $1\leq i\leq m$, we
obtain $\P_\theta (\FDP_m(\tau^0_{\hat{{\ell}}}) >\alpha
)>\zeta+
10^{-3}$. Admittedly, the FDP control is just slightly violated.
Nevertheless, this gives numerical evidence of the following fact.
%

\begin{fact}\label{factoracleSD}
RW's heuristic does not always provide a valid FDP control for finite
$m$ in its oracle step-down form in model \eqref{equi-correlated}.
\end{fact}

Note that the case of the \textit{non-oracle} adaptive version is studied with extensive
simulations in Section~S-6 in \hyperref[suppA]{Supplementary Material}, the conclusion is similar.
Fact~\ref{factoracleSD} is interesting from a theoretical point of
view: it annihilates any hope of finding a general finite sample proof
of FDP control in the step-down case, even under a very simple form of
positive dependence.

\subsection{Existing modifications}\label{existingmodif}

Facts~\ref{factoracleSU} and~\ref{factoracleSD} indicate that to obtain
a provable finite sample control, it is appropriate to slightly
decrease the initial $k$-FWE-based critical values $(\tau_{\ell
})_{1\leq
{\ell}
\leq m}$.
Interestingly, several existing procedures that provably control the
FDP can be reinterpreted as modifications of the $\tau_{\ell}$'s.
In the literature, we have identified the following principles
that provide a control of the FDP under general dependence:
\begin{itemize}
\item The ``diminution'' principle \citet{GHS2013},
\citeauthor{RS2006} (\citeyear{RS2006,RS2006b}): first,
establish a rigorous upper-bound for $\P(\FDP>\alpha)$ for a step-down
or step-up procedure with arbitrary critical values $(c_{\ell
}(x))_{1\leq
{\ell}\leq m}$ depending on a single parameter $x$. Second, adjust $x$ to
make the bound smaller than $\zeta$. As an illustration,
\citeauthor{RS2006} (\citeyear{RS2006,RS2006b}) have proposed the following bound that can be
rewritten as follows (see Section~\ref{sec:newbound} for a proof):
%
%
\begin{equation}
C^{\mathrm{RS}}(x)= \max_{1\leq u\leq m} \Biggl\{u \sum
_{{\ell}=1}^{b_\alpha(u)} \frac
{c_{\ell}(x)-c_{{\ell}-1}(x)}{d({\ell},m,u)} \Biggr\} ;\label{equCRS}
\end{equation}
where for all $u,{\ell}\in\{0,\ldots,m\}$, we let
%
%
\begin{eqnarray}
\label{balphau} b_\alpha(u) &=&\cases{ %
\bigl(\bigl\lfloor(m-u)/(1-\alpha)\bigr\rfloor+1\bigr)\wedge\bigl(\lceil u/\alpha
\rceil-1\bigr)\wedge m,\vspace*{2pt}\cr
\quad\hspace*{91pt}\mbox{(step-down)},
\vspace*{2pt}\cr
{}\bigl(\lceil u/\alpha\rceil -1\bigr)\wedge m,\qquad\mbox{(step-up)},}
\\[3pt]
d({\ell},m,u) &=&\cases{ %
\lfloor\alpha{\ell}
\rfloor+1,& \quad$\mbox{(step-down)},$
\vspace*{2pt}\cr
\bigl(\lfloor\alpha{\ell}\rfloor+1\bigr)\vee({\ell}-m+u),&\quad
$\mbox{(step-up)}.$}\label{dlmu}
\end{eqnarray}
This bound does not incorporate the dependence. Finally, let us mention
that the diminution principle has been recently followed by using much
more sophisticated bounds that incorporate the pairwise dependence; see
Theorems~3.7 and~3.8 in \citet{GHS2013}.
\item The ``augmentation'' principle \citet{DLP2004b,Far2009}: consider
the $1$-FWE controlling procedure at level $\zeta$ rejecting the null
hypotheses corresponding to the set
$
\mtc{R}^{(1)}=\{1\leq i\leq m \dvtx p_i\leq\tau_1(\zeta)\}
$, denote ${\ell}^{(1)}$ the number of rejections of $\mtc{R}^{(1)}$ and
\[
\widetilde{{\ell}}^{\mathrm{Aug}} = \bigl\lfloor{\ell}^{(1)}/(1-\alpha)
\bigr\rfloor \wedge m.
\]
Then the ``augmented'' procedure rejects the nulls associated to the
$\widetilde{{\ell}}^{\mathrm{Aug}}$ smallest $p$-values.
This procedure can incorporate the dependence if $\mtc{R}^{(1)}$ is
appropriately chosen.

\item The ``simultaneous'' $k$-FWE control proposed in \citet{GW2006}:
consider critical values $(\tau_{{\ell}}(\zeta/m))_{1\leq{\ell
}\leq m}$ (with
$\zeta$ divided by $m$), and let
\[
\widetilde{{\ell}}^{\mathrm{sim}} = \biggl\lfloor\frac{\max \{ R(\tau
_{{\ell}
}(\zeta
/m)) - \lfloor\alpha{\ell}\rfloor\dvtx{\ell}\leq R(\tau_{{\ell
}}(\zeta
/m)), {\ell}
\geq0 \}}{1-\alpha} \biggr
\rfloor\wedge m.
\]
Then the ``simultaneous'' procedure rejects the nulls corresponding to
the $\widetilde{{\ell}}^{\mathrm{sim}}$ smallest $p$-values. Again, this procedure
is able to incorporate the dependence if the $\tau_{{\ell}}$'s are
suitably built.
\end{itemize}

\subsection{Two new modifications}\label{secnewmodif}

This section presents new results that can be seen as modifications of
$k$-FWE based procedures that ensure finite sample FDP control. Both
modifications incorporate the dependence between the $p$-values.
Furthermore, the numerical experiments of Section~\ref{sec:num} show
that they are more powerful than the state-of-the-art procedures
described in Section~\ref{existingmodif}.

\textit{A first modification}. The first result follows the
``diminution'' principle. For any arbitrary critical values $(c_{\ell}
(x))_{1\leq{\ell}\leq m}$ (depending on a variable $x$), let
$C^{\mathrm{ex}}(x)$ be
%
%
\begin{eqnarray}\label{equC}
&& \max_{1\leq u \leq m} \Biggl\{\sum_{{\ell}=1}^{b_\alpha(u)}
\mathop{\max_{\theta\in\Theta}}_{ m_0(\theta)=u} \bigl\{ \bigl(\P
_{\theta} \bigl( V_m\bigl(c_{\ell}(x) \bigr)\geq d({
\ell}-1,m,u) \bigr)\nonumber\\
&&\hspace*{95pt}{}- \P _{\theta
} \bigl( V_m
\bigl(c_{{\ell}-1}(x)\bigr)\geq d({\ell}-1,m,u) \bigr) \bigr)
\\
&&{}\hspace*{93pt}  \wedge \bigl(\P_{\theta} \bigl( V_m
\bigl(c_{\ell}(x)\bigr)\geq d({\ell} ,m,u) \bigr)\nonumber\\
&&\hspace*{113pt}{}- \P_{\theta}
\bigl( V_m\bigl(c_{{\ell}-1}(x)\bigr)\geq d({\ell},m,u) \bigr)
\bigr) \bigr\} \Biggr\},\nonumber
\end{eqnarray}
where $b_\alpha(u)$ and $d({\ell},m,u)$ are given by \eqref{balphau} and
\eqref{dlmu}, respectively. The following result is established in
Section~\ref{sec:newbound}.

%
\begin{theorem}\label{main-th1}
Let us consider either the fixed model ($\Theta=\Theta^F$) or the
uniform model ($\Theta=\Theta^U$) and
any family of critical values\vspace*{1pt} $(c_{\ell}(x))_{1\leq{\ell}\leq m}$,
$x\geq0$,
such that $c_m(0)=0$.
Consider some $x^\star\geq0$ satisfying $C^{\mathrm{\textup{ex}}}(x^\star)\wedge
C^{\mathrm{RS}}(x^\star)\leq\zeta$, where $C^{\mathrm{ex}}(\cdot)$ is defined by
\eqref
{equC} and $C^{\mathrm{RS}}(\cdot)$ by \eqref{equCRS}. Let $\wh{{\ell}}$ be the
number of rejections of the step-down \eqref{SD-classic} [resp., step-up
\eqref{SU-classic}] algorithm associated to the critical values
$(c_{\ell}
(x^\star))_{1\leq{\ell}\leq m}$.
Then the FDP control \eqref{FDPcontrol} holds, with $\wh{t}_{m}=\tau
_{\hat{\ell}}$.
\end{theorem}

Theorem~\ref{main-th1} can be applied with any starting critical values
$(c_{\ell}(x))_{1\leq{\ell}\leq m}$. A~choice in accordance with RW's
heuristic is
$c_{\ell}(x) = x \wt{\tau}_{\ell}$, ${\ell}\in\{1,\ldots,m\}$,
$x\geq0$,
where $(\wt{\tau}_{\ell})_{1\leq{\ell}\leq m}$ are the adaptive
$k$-FWE based
critical values \eqref{equtaukadapt} for some appropriate bounding device.
Next, while Theorem~\ref{main-th1} does not require any assumption on
the dependence, it implicitly assumes that the function $C^{\mathrm{ex}}(\cdot)$
is known or easily computable. This is the case, for instance, in the
model \eqref{facmodel} because we have
%
%
\begin{eqnarray}\label{equCrhoExact}
C^{\mathrm{ex}}(x)& =& \max_{1\leq u \leq m} \Biggl\{\sum
_{{\ell}=1}^{b_\alpha(u)} \bigl( {B}^0_m
\bigl(c_{\ell}(x) ,d({\ell}-1,m,u),u\bigr)\nonumber\\
&&\hspace*{58pt}{} -{B}^0_m
\bigl(c_{{\ell
}-1}(x) ,d({\ell} -1,m,u),u\bigr) \bigr)
\nonumber
\\[-8pt]
\\[-8pt]
\nonumber
&&{}\hspace*{56pt}\wedge \bigl( {B}^0_m\bigl(c_{\ell}(x) ,d({
\ell},m,u),u\bigr)\\
&&\hspace*{79pt}{} -{B}^0_m\bigl(c_{{\ell}-1}(x) ,d({
\ell} ,m,u),u\bigr) \bigr) \Biggr\} ,\nonumber
\end{eqnarray}
where ${B}^0_m(t,k,u)$ is the exact bounding device defined by \eqref
{equ-multirhoExact}.
A second illustration is the Gaussian case where $\Gamma$ is known but
arbitrary and where the model is $\Theta=\Theta^U$. In this situation,
$C^{\mathrm{ex}}(x)$ in \eqref{equC} can be approximated by Monte Carlo calculations.
Finally, let us underline that Theorem~\ref{main-th1} provides FDP
control even if the incorporated dependence is not positive.

\textit{A second modification}. The second result presented in this
section relies on
the $K$-Markov device $B^0_m(t,k,u)$ given by \eqref{KMarkovbounding}
(for some integer $K\geq1$).\vadjust{\goodbreak} It specifically uses the two following
assumptions (here $\Theta=\Theta^F$ only):
%
%
\renewcommand{\theequation}{Exch-${\mathcal{H}}_0$}
\begin{equation}\label{exchangeability}\qquad\quad
\begin{tabular}{p{280pt}@{}}
for all $\theta\in\Theta$, for any permutation
$\sigma$ of $\{1,\ldots,m\}$ with $\sigma(i)=i$ for all $i\notin {\mtc{H}}
_0(\theta)$,
the distribution of $(p_{\sigma(i)})_{1\leq i
\leq m}$ is equal to the
one of $(p_i)_{1\leq i \leq m}$;
\end{tabular}\vspace*{-12pt}
\end{equation}
\renewcommand{\theequation}{Posdep}
\begin{equation}\label{weakPRD}\qquad\quad
\begin{tabular}{p{280pt}@{}}
for all $\theta\in\Theta$, for any
measurable nondecreasing set $D\subset[0,1]^m$ and subset $X\subset{
\mtc {H}}_0(\theta )$,
$x\in[0,1]\mapsto\P_\theta \bigl((p_i)_{1\leq i \leq
m}
\in D\cond\forall i\in X, p_i\leq x \bigr)$ is nondecreasing.
\end{tabular}
\end{equation}
In \eqref{weakPRD}, a set $D\subset[0,1]^m$ is said nondecreasing if
for any $x,y\in[0,1]^m$ such that $x\in D$, the inequality $x\leq y $
(holding component-wise) entails $y\in D$. Condition \eqref{weakPRD}
induces a form of positive dependence between the $p$-values. It is
stronger than the condition of positive dependence ensuring FDR control
for the BH procedure, for which the conditioning holds w.r.t. only one
element; see \citet{BY2001,BR2008,BDR2014}.
However,
assumption \eqref{weakPRD} is satisfied as soon as the $p$-value
family is
multivariate totally positive of order 2 (MTP2); see \citet{Sar1969}.
We refer to Section~3 of \citet{KR1981} for several examples of MTP2
models, which thus satisfy assumption~\eqref{weakPRD}. More explicit
examples will be provided at the end of the section.

Now, assuming \eqref{exchangeability}, we consider for ${\ell}\in\{
1,\ldots
,m\}$,
%
%
\renewcommand{\theequation}{\arabic{equation}}
\setcounter{equation}{25}
\begin{equation}
\label{newcritval} {\tau}^{\mathrm{new}}_{\ell}=\cases{ %
 \wt{\tau}_{\ell}\bigl( \lambda\zeta,m({\ell})\bigr), &\quad
$\mbox{if ${\ell }\geq{\ell}_K$},$
\vspace*{2pt}\cr
\displaystyle\biggl(\frac{(1-\lambda)\zeta(\lfloor\alpha{\ell}\rfloor+1)
}{m({\ell})
} \biggr) \wedge\wt{\tau}_{{\ell}_K}( \lambda
\zeta,m), &\quad $\mbox{if ${\ell }< {\ell} _K$,}$}
\end{equation}
%
where $\l_K=\lceil (K-1)/\alpha\rceil$, $\lambda\in [0,1]$ is some tuning parameter.
Also, $\wt{\tau}_\l(\zeta,u)$ denotes the value of $t$ obtained by solving the equation
%
%
\begin{eqnarray}\label{equ-Kmarkov}
&&\mathop{\sup_{\theta\in\Theta}}_{ m_0(\theta)=u} \Bigl\{ P_\theta
\Bigl(\max_{ i\in X_0}\{ p_i\}\leq t \Bigr) \Bigr\}
\nonumber
\\[-8pt]
\\[-8pt]
\nonumber
&&\qquad=
\zeta\frac
{(\lfloor
\alpha{\ell}\rfloor+1)(\lfloor\alpha{\ell}\rfloor+1-1)\cdots
(\lfloor
\alpha{\ell}\rfloor+1-K+1)} {u(u-1)\cdots(u-K+1)},
\end{eqnarray}
where $X_0$ denotes any subset of ${\mtc{H}}_0$ of cardinal $K$. The
following result holds; see Section~\ref{proofth-newcontrol} for a proof.

%
\begin{theorem}\label{th-newcontrol}
In the fixed model $\Theta=\Theta^F$, let $\wh{{\ell}}$ be the
number of
rejections of the step-up \eqref{SU-classic} algorithm associated to
the critical values\break $({\tau}^{\mathrm{new}}_{\ell})_{1\leq{\ell}\leq m}$
given by
\eqref{newcritval}. Then the finite sample FDP control \eqref
{FDPcontrol} holds for $\wh{t}_{m}=\tau^{\mathrm{new}}_{\hat{{\ell}}}$ under
assumptions \eqref{weakPRD} and \eqref{exchangeability}.
\end{theorem}

The proof of Theorem~\ref{th-newcontrol} is given in Section~\ref
{proofth-newcontrol}. It shares some similarity with the proofs
developed in \citet{Sar2007} in the FDR case.
When $K=1$ and $\lambda=1$, the critical values $({\tau}^{\mathrm{new}}_{\ell}
)_{1\leq{\ell}\leq m}$ are the Lehmann--Romano critical values \eqref
{criticalvaluesMarkov}, and thus
Theorem~\ref{th-newcontrol} is in accordance with Proposition~\ref
{RWalreadyproved}(ii) and Theorem~3.1 of \citet{GHS2013} because
Simes's inequality is valid in that case.
The originality of Theorem~\ref{th-newcontrol} lies in the case $K>1$
that allows us to incorporate the dependence in an FDP controlling
procedure. Below, some examples are provided in the one-sided location
models \eqref{equonesided}:
\begin{longlist}[(ii)]
\item[(i)] When the noise $Y$ is Gaussian multivariate,
assumption~\eqref{exchangeability} imposes equicorrelation between the
$p$-values $(p_i,i\in{\mtc{H}}_0)$, say $\rho$-\break equicorrelation with
$\rho
\in
[0,1)$. In this case, equation~\eqref{equ-Kmarkov} can be solved by using
$
\P_\theta ( \forall i\in X_0, p_i\leq t ) = \E[
(F_{0}(t,W))^{K}]$,
where $F_{0}$ is given by\break \eqref{Fzero-equi-correlated} and $W\sim
\mtc{N}(0,1)$.
Furthermore, the $p$-value family is\break MTP2 if and only if $-\Gamma^{-1}$
has nonnegative off-diagonal elements; see, for example, \citet{Rin2006}.
For instance, both assumptions are satisfied if $\Gamma$ is $\rho
$-equi-correlated \eqref{equi-correlated}. Additional examples can be
provided with matrices $\Gamma$ such that $(\Gamma_{i,j})_{i\in{\mtc{H}}
_0,j\in{\mtc{H}}_0}$ is $\rho$-equi-correlated while $-\Gamma^{-1}$ has
nonnegative off-diagonal elements.

\item[(ii)] Consider \eqref{facmodel} in the particular case where
$X_i=\mu_i+ c_i W + \zeta_i- (a/b) w_0$,
where $c_1\sim\gamma(a,b)$, $W$ is a positive random variable
($w_0=\E W$) and $\zeta_1$ is centered with a log-concave density.
In this case, the $p$-value family is MTP2 by Proposition~3.7 and~3.9
of \citet{KR1980}, which entails \eqref{weakPRD}. Assumption~\eqref
{exchangeability} also clearly holds, and the LHS of \eqref
{equ-Kmarkov} is
$
\E[ (F_{0}(t,W))^{K}]$,
where $F_0$ is defined by \eqref{Fzero-facmodel}.
\end{longlist}

%
\begin{remark}
Assumption~\eqref{weakPRD} is, strictly speaking, weaker than MTP2
property. For instance, \eqref{weakPRD} is satisfied in the Gaussian
case where $\Gamma_{i,j}\geq0$ for $i\in{\mtc{H}}_1$ and $j\in
{\mtc{H}}_0$ and
$\Gamma_{i,j}=1$ for $i,j\in{\mtc{H}}_0$.
\end{remark}

\section{Asymptotic results}\label{sec:asymp}

The goal of this section is to study RW's heuristic from an asymptotic
point of view.

\subsection{Setting and assumptions}\label{settingasymp}

In this section, the FDP control under study is asymptotic: we search
$\wh{t}_m$ such that
%
%
\begin{equation}
\limsup_m \bigl\{\P_{\theta^{(m)}} \bigl(
\FDP_m (\wh {t}_m )> \alpha \bigr) \bigr\}\leq
\zeta\label{AsympsFDPcontrol}.
\end{equation}
This requires us to consider a \textit{sequence} of models $(\Theta
^{(m)}, m\geq1)$ (fixed mixture models here) and a sequence of
parameters $(\theta^{(m)}, m\geq1)$ with $\theta^{(m)}\in\Theta
^{(m)}$ for all $m\geq1$. The latter sequence is assumed to be fixed
once for all throughout this section.
Moreover, we will assume throughout this study the following common assumption:
%
%
\begin{equation}
m_0\bigl(\theta^{(m)}\bigr)/m \rightarrow\pi_0\qquad
\mbox{where } \pi_0\in(0,1). \label{condm0}
\end{equation}
In particular, any sparse situation where $m_0(\theta
^{(m)})/m\rightarrow1$ is excluded. Also, under~\eqref{condm0}, we let
$\pi_1=1-\pi_0\in(0,1)$.

Useful assumptions on $(\theta^{(m)}, m\geq1)$ are the following weak
dependence assumptions on the processes ${\widehat{\mathbb
{G}}}_{m}(t)=R_m(t)/m$ and
${\widehat{\mathbb{G}}}_{0,m}
(t)=V_m(t)/m_0$, $t\in[0,1]$:
%
%
\renewcommand{\theequation}{weakdep}
\begin{equation}\qquad\quad 
\hspace*{24pt} \Vert{\widehat{\mathbb{G}}}_{m}- G\Vert_\infty=
o_P(1)\hspace*{20pt} \mbox{for some continuous $G\dvtx[0,1]\to[0,1]$; }
\label{weakdep}\vspace*{-12pt}
\end{equation}
\renewcommand{\theequation}{weakdep0}
\begin{equation}
\hspace*{-22pt}\Vert{\widehat{\mathbb{G}}}_{0,m}- I
\Vert_\infty= o_P(1)\qquad \mbox {for $I(t)=t$, $t\in [0,1]$. }
\label{weakdep0} 
\end{equation}
These weak dependence conditions are widely used in the context of
multiple testing; see, for example, \citet{STS2004,GF2013} and the
stronger condition \eqref{FCLT} further on.
In the particular one-sided Gaussian multivariate setting, these
conditions have been studied in \citet{SL2011,FHG2012,DR2014} (among others).
Lemma S-3.1 
of the\break \hyperref[suppA]{Supplementary Material} states that
assumptions~\eqref{weakdep} and \eqref{weakdep0} are satisfied with
$G(t)=\pi_0 t + \pi_1 F_1(t)$ and $F_1(t)=\int_0^\infty\ol{\Phi}(
{\ol{\Phi}}^{-1}
(t)-\beta) \,d\nu(\beta)$ under \eqref{condm0} if the following
conditions hold:
%
%
\renewcommand{\theequation}{\mbox{Conv-alt}}
\begin{eqnarray}\label{Cmu}
(m_1)^{-1} \sum_{i=1}^m
H_i \delta_{\mu_i} \mathop{\longrightarrow }^{\mathrm{weak}}
\nu
\nonumber
\\[-8pt]
\\[-8pt]
\eqntext{\mbox{for a distribution $\nu$ on
$\R^+$ with $\nu \bigl(\{0\} \bigr)=0$, }}\vspace*{-12pt}
\end{eqnarray}
\renewcommand{\theequation}{weakdepGauss}
\begin{equation}
m^{-2}\sum_{i,j=1}^m
(\Gamma_{i,j})^2 \rightarrow0.\label
{weakdepGauss}
\end{equation}
Also, let us underline that under \eqref{equi-correlated},
assumption~\eqref{weakdepGauss} is satisfied whenever $\rho=\rho
_m\to0$.

Finally, we also explore in this section \textit{strong dependence},
through the factor model \eqref{facmodel}. This includes \eqref
{equi-correlated} for a parameter $\rho\in(0,1)$ taken fixed with $m$.

\subsection{The BH procedure and FDP control}

Let us go back to Figure~\ref{FDPdep}.
When $\rho=0$, even if the BH procedure is only intended to control the
expectation of the FDP at level $\alpha$, the $95\%$ quantile of the
FDP is still close to $\alpha$. This comes from the concentration of
the FDP of the BH procedure around $\pi_0 \alpha<\alpha=0.2$ as $m$
grows to infinity. It is well known that this quantile converges to
$\pi
_0 \alpha$ as $m$ grows to infinity, so that the limit in \eqref
{AsympsFDPcontrol} is equal to zero; see, for example, \citet{Neu2008}.
In other words, the FDP concentration combined with the slight amount
of conservativeness due to $\pi_0<1$ ``prevents'' the FDP from exceeding
$\alpha$.
The consequence is simple: the BH procedure controls the FDP
asymptotically in the sense of \eqref{AsympsFDPcontrol} under
independence. As a matter of fact, the latter also holds under weak dependence.

%
\begin{lemma}\label{BHFDP}
Consider the BH procedure, that is, the step-up procedure \eqref
{SU-classic} associated to the linear critical values $\tau_{\ell
}=\alpha
{\ell}
/m$, $1\leq{\ell}\leq m$.
Assume that $(\theta^{(m)}, m\geq1)$ satisfies \eqref{condm0},
\eqref
{weakdep}, \eqref{weakdep0} and further assume that $G$ satisfies the
following property:
%
%
\renewcommand{\theequation}{Exists}
\begin{equation}
\label{exists} 
\mbox{ there exists } t\in(0,1) ,\mbox{ such that }
G(t)>t/\alpha.
\end{equation}
Then we have $\P_{\theta^{(m)}} (\FDP_m (\tau_{\hat{\ell}
} )>
\alpha )\rightarrow0$.
\end{lemma}

Although this result seems new, its proof, provided in Section~\ref
{sec:proofBHFDP}, can certainly be considered as standard; see, for
example, \citet{GW2004,FDR2007}.
Also, while Lemma~\ref{BHFDP} does not require any Gaussian assumption
in general, all the assumptions of Lemma~\ref{BHFDP} are satisfied
under \eqref{condm0}, \eqref{Cmu} and \eqref{weakdepGauss}.

In the literature, even under independence, it is common to exclude the
BH procedure while studying \eqref{AsympsFDPcontrol}.
For instance, Proposition~4.1 in \citet{CT2008} shows that the ``oracle''
version of the BH procedure, that is, the step-up procedure with
critical values $\alpha{\ell}/m_0$, ${\ell}\in\{1,\ldots,m\}$, has
a FDP
exceeding $\alpha$ with a probability tending to $1/2$. Since the
oracle BH procedure is often considered to be better than the original
BH procedure, it is thus tempting to exclude the BH procedure when
studying an FDP control of type \eqref{AsympsFDPcontrol}.
Lemma~\ref{BHFDP} shows that, perhaps surprisingly, this is a mistake:
BH procedure is interesting when providing \eqref{AsympsFDPcontrol} and
does not suffer from the same drawback as the oracle BH procedure.

By contrast, if the dependence is not weak, the BH procedure can fail
to control the FDP as $m$ is tending to infinity. Under \eqref
{equi-correlated} and when the $p$-values under the alternative are
zero (Dirac uniform configuration), this fact has been formally
established in Theorem~2.1 of \citet{FDR2007}, by showing that the limit
of the FDP of the BH procedure is not deterministic anymore and hence
can exceed $\alpha$ with a positive probability, which is obviously not
related to $\zeta$ (because BH critical values do not depend on $\zeta$).

\subsection{RW's heuristic under weak dependence}

The two results provided in this section both validate the use of RW's
heuristic under weak dependence.
Since the BH procedure is valid in this case, they are of limited
interest in practice. Nevertheless, we believe that they suitably
complement our overview on RW's heuristic.
The first result is proved in Section~\ref{proof:thRwheuriweakdep0} via
technics similar to those used for proving Lemma~\ref{BHFDP}.

%
\begin{theorem}\label{thRwheuriweakdep0}
Consider the one-sided location model \eqref{equonesided} with the
full null process $V_m'(\cdot)$ being defined by \eqref{fullnullprocess}.
Assume that the (nonadaptive) exact bounding device is such that for
all $t,k$ $\ol{B}_m(t,k)=\P_\theta(V_m'(t)\geq k)$ for the parameter
$\theta=\theta^{(m)}$ at hand.
Consider the critical values $\ol{\tau}_{\ell}$, $1\leq{\ell}\leq m$,
derived from $\ol{B}_m$ as in \eqref{equtauknonadapt}, and consider the
corresponding step-up procedure \eqref{SU-classic}.
Assume that $(\theta^{(m)}, m\geq1)$ satisfies \eqref{condm0},
\eqref
{weakdep}, \eqref{weakdep0}, \eqref{exists} and that $V'_m(t)/m$
converges in probability to $t$ for any $t\in[0,1]$ (i.e., weak
dependence for the full null process).
Then we have $\P_{\theta^{(m)}} (\FDP_m (\ol{\tau}_{\hat
{\ell}
} )> \alpha )\rightarrow0$.
\end{theorem}

The above result shows that RW's procedure used with the exact bounding
device turns out to have an asymptotic exceedance probability of zero
under weak dependence, likewise the BH procedure. Again, this is due to
the convergence of the FDP toward $\pi_0\alpha<\alpha$.
Hence, perhaps disappointingly, $\zeta$ plays no role in the limit,
which indicates that using the simpler BH procedure seems more
appropriate in this case.

Nevertheless, when $m_0$ is known, an interesting point is that, while
the oracle BH procedure fails to control the FDP (as discussed in the
above section), oracle RW's method maintains the FDP control.
To show this, we need to slightly strengthen the assumption \eqref
{weakdep0} by assuming
that the following central limit theorem holds for the (rescaled)
process $V_m(\cdot)$:
%
%
\renewcommand{\theequation}{FLT}
\begin{equation}\qquad
\left\{ %
\begin{minipage} {10cm} There is a rate $r_m
\rightarrow\infty$ such that the process $Z_m(t)=r_m
\bigl(V_m(t)/m-\bigl(m_0(m)/m\bigr) t \bigr)$ satisfies,
for any $K=[a,b]\subset(0,1)$, the convergence $\bigl(Z_m(t)
\bigr)_{t\in K} \leadsto\bigl(Z(t)\bigr)_{t\in K}$ (for the
Skorokhod topology), for a process $\bigl(Z(t)\bigr)_{t\in K}$ with
continuous paths and such that the random variable $Z(t)$ has a continuous
increasing c.d.f. for all $t\in K$. \end{minipage} \right.
\label{FCLT} 
\end{equation}

For instance, under \eqref{condm0}, assumption~\eqref{FCLT} holds when
the $p$-values $(p_i, H_i=0)$ are i.i.d. by Donsker's theorem. More
generally, dependencies satisfying ``mixing'' conditions can also lead
to \eqref{FCLT}; see, for example, \citet{DP2007,DLST2010} or \citet
{Far2007}. Recently, some efforts have been undertaken to consider
other types of dependence, not necessarily locally structured; see
\citet
{Sou2001,BS2013}. In the case of a Gaussian multivariate structure,
explicit sufficient conditions on $\Gamma$ are provided in \citet
{DR2014}. The following result is proved in Section~\ref{sec:thRwheuriweakdep}.

%
\begin{theorem}\label{thRwheuriweakdep}
Assume that the exact oracle bounding device is such that for all $t,k$
$B_m^0(t,k,m_0)=\P_{\theta}(V_m(t)\geq k)$ for the parameter $\theta
=\theta^{(m)}$ at hand.
Consider the critical values $\tau^0_{\ell}$, $1\leq{\ell}\leq m$, derived
from $B_m^0$ as in \eqref{equtaukoracle}, and consider the
corresponding step-up procedure \eqref{SU-classic}. Assume that
$(\theta
^{(m)}, m\geq1)$ satisfies \eqref{condm0} with $\pi_0>\alpha$,
\eqref
{weakdep}, \eqref{FCLT} and further assume that $G$ satisfies the
following property:
%
%
\renewcommand{\theequation}{Unique}
\begin{equation}
G(t)=\pi_0 t/\alpha\mbox{ has at most one solution on $(0,1)$;}
\label
{unique}\vspace*{-12pt}
\end{equation}
\renewcommand{\theequation}{NonCritical}
\begin{equation}
\lim_{t\to0^+} G(t)/t \in(\pi_0/\alpha,+
\infty]. \label
{noncritical} 
\end{equation}
Then we have $\P_{\theta^{(m)}} (\FDP_m ({\tau}^0_{\hat
{\ell}
} )> \alpha )\rightarrow\zeta$.
\end{theorem}

Roughly speaking, the essence of the argumentation is as follows:
when $\wh{t}$ converges in probability to some deterministic quantities,
then the fluctuations of $\wh{\l}/m$
asymptotically disappear in probability \eqref{heuristic}, and thus the
latter is equal to $\zeta$ by definition of the oracle exact bounding device.
Note that a similar reasoning has been made at the end of Section~7 in
\citet{GW2006}.
Here, we derive sufficient conditions that make this informal argument rigorous.

Markedly, in Theorem~\ref{thRwheuriweakdep}, the limit of the
probability is exactly $\zeta$; hence there is no loss in the level of
RW's method. However, since $m_0$ is often unknown (and seems hard to
estimate at a rate faster than $r_m$), the interest of this result
remains mainly theoretical.
Finally note that \eqref{unique} and \eqref{noncritical} are classical
conditions when studying asymptotic properties of step-up procedures;
see, for example, \citet{GW2002,Chi2007,Neu2008}.

\subsection{RW's heuristic under strong dependence}\label{sec:strongdep}

Here, we study the asymptotic properties of RW's method under strong dependence
by focusing on models of the type \eqref{facmodel}.
A crucial assumption is as follows:
%
%
\renewcommand{\theequation}{Posdep-facmod}
\begin{equation}
\label{asympposdep} 
\qquad\qquad\quad\qquad\mbox{for any $t\in(0,1)$, the function $w\mapsto
F_{0}(t,w)$ is increasing,}
\end{equation}
where $F_0$ is given by \eqref{Fzero-facmodel}. The latter assumption
is a form of positive dependence which is specific to \eqref{facmodel}:
it roughly means that the variable $W$ disturbs each $p$-value
distribution in an ``unidirectional'' manner. For instance,
\eqref{asympposdep} is satisfied if $\P(c_1\geq0)=1$, $\P(c_1>0)>0$
and $\xi_1$ has a distribution function which is continuous increasing
on $\R$.
As an illustration, \eqref{asympposdep} is satisfied under \eqref
{equi-correlated} ($\rho> 0$) but not necessarily under \eqref
{alt-equi-correlated}.

\textit{Asymptotic view of RW's heuristic.}
Under appropriate assumptions, the exact device \eqref
{equ-multirhoExact} is such that, for a sequence $k_m$ with
$k_m/m\rightarrow\kappa$,
%
%
\renewcommand{\theequation}{\arabic{equation}}
\setcounter{equation}{29}
\begin{equation}
\ol{B}_m(t,k_m) = \P( N_m/m \geq
k_m/m) \rightarrow\P_W \bigl(F_{0}(t,W)\geq
\kappa\bigr),\label{equasympRWheu}
\end{equation}
where $F_{0}$ is defined by \eqref{Fzero-facmodel}, and $N_m$ follows a
binomial distribution of parameters $m$ and $F_{0}(t,W)$, conditionally
on $W$.
By taking $\kappa=F_{0}(t,q_\zeta)$ where $q_\zeta$ is such that $\P
(W\geq q_\zeta)\leq\zeta$,
the probability on the RHS of \eqref{equasympRWheu} is smaller than or
equal to
$\P_W (F_{0}(t,W)\geq F_{0}(t,q_\zeta) ) \leq\P_W (W\geq q_\zeta)
\leq
\zeta$, provided that \eqref{asympposdep} holds.
Now, RW's heuristic (taken in an asymptotic sense)
leads to the following equation for the critical values:
%
%
\begin{equation}
\label{equ-critvalues-RWasymp} F_{0}(\tau_{\ell},q_\zeta) = \alpha{
\ell}/m,\qquad 1\leq{\ell}\leq m.
\end{equation}
Under independence, this gives the BH critical values. In the
equi-correlated case \eqref{equi-correlated} (with $\rho\geq0$), this
yields the critical values \eqref{simpletaul} mentioned in the
\hyperref[sec1]{Introduction} of the paper.

\textit{A modification based on DKW's concentration inequality.}
The following result shows that a simple modification of \eqref
{equ-critvalues-RWasymp} provides FDP control (see Section~\ref
{sec:proofwithDKW} for a proof).

%
\begin{theorem}\label{thmodifasymp}
Let $\lambda\in(0,1)$. In a model \eqref{facmodel} satisfying\break \eqref
{asympposdep}, consider any critical values $\tau_{\ell}$, $1\leq
{\ell}\leq
m$, satisfying
%
%
\begin{equation}
\label{equ-critvalues-RWasymp-modif} F_{0}(\tau_{\ell},q_{\zeta(1-\lambda)}) \leq
\biggl( \frac{\alpha
{\ell}
}{m} - \biggl( \frac{-\log(\lambda\zeta/2)}{2m} \biggr)^{1/2}
\biggr)_+,\qquad 1\leq {\ell}\leq m,
\end{equation}
where for any $x\in(0,1)$, $q_x\in\R$ is such that $\P(W\geq
q_x)\leq
x$ and $F_{0}$ is defined by \eqref{Fzero-facmodel}.
Consider the step-up procedure \eqref{SU-classic} associated to the
critical values $\tau_{\ell}$, ${\ell}=1,\ldots,m$. Then it
controls the FDP;
that is, \eqref{FDPcontrol} holds with $\wh{t}=\tau_{\hat{\ell}}$.
\end{theorem}

While this result is nonasymptotic, it is intended to be used for
large values of $m$ in order to reduce the influence of the remainder terms.
However, even for large values of $m$, \eqref
{equ-critvalues-RWasymp-modif} imposes to set the first critical values
to zero, which may be undesirable. The next section presents conditions
allowing to drop these annoying remainder terms.

\textit{An asymptotic validation of RW's heuristic.}
Here, we present situations for which the raw critical values \eqref
{equ-critvalues-RWasymp} can be used to get an asymptotic FDP control.
We consider the following additional distribution assumptions on $c_1$,
$W$ and $\xi_1$ in \eqref{facmodel}:
\begin{longlist}[(iii)]
\item[(i)] $c_1$ is a random variable with a finite support in $\R^+$;
\item[(ii)] the distribution function of $W$ is continuous;
\item[(iii)] the function $x\in\R\mapsto\ol{F}_{\xi}(x)=\P(\xi
_1\geq
x)$ is continuous increasing and is such that, for all $y\in\R$, as
$x\rightarrow+\infty$,
%
%
\begin{equation}
\label{LRseparate} \frac{\ol{F}_{\xi}(x-y) }{ \ol{F}_\xi(x) }\rightarrow
\cases{ +\infty,&\quad $\mbox{if $y>0$},$
\vspace*{2pt}\cr
0, &\quad $\mbox{if $y<0.$}$}
\end{equation}
\end{longlist}
When $\xi_1$ has a log-concave density, condition \eqref{LRseparate}
can be reformulated in terms of a density ratio; see, for example, the
relations in Section~S-5 of \citet{NR2012}.
For instance, a simple class of distributions satisfying \eqref
{LRseparate} are the so-called \textit{Subbotin} distributions, for which
the density of $\xi_1$ is given by $e^{-|x|^\gamma/\gamma}$ (up to a
constant), for some parameter $\gamma>1$; see Section~5 of \citet
{Neu2013} for more details on this.

The following result is proved in Section~\ref{sec:proofthasympfact}.

%
\begin{theorem}\label{thasympfact}
Consider the one-sided testing problem \eqref{equonesided} with all
alternative means equal to some $\beta>0$, and assume that $\theta
^{(m)}$ satisfies \eqref{condm0}. Consider \eqref{facmodel} satisfying
the assumptions \textup{(i)}, \textup{(ii)} and \textup{(iii)} above.
Consider the step-up procedure \eqref{SU-classic} associated to the
critical values $\tau_{\ell}$, ${\ell}=1,\ldots,m$ satisfying
\eqref
{equ-critvalues-RWasymp}.
Then the asymptotic FDP control \eqref{AsympsFDPcontrol} holds with
$\wh
{t}_m=\tau_{\hat{\ell}}$.
\end{theorem}

%
\begin{table}[b]
\caption{Procedures used in Figures~\protect\ref{powernotusingrho}
and~\protect\ref{powerusingrho}; see Sections~\protect\ref
{existingmodif} and~\protect\ref{secnewmodif} for more
details. All the procedures are step-up; ``e.b.d.'' means ``exact
bounding device''}\label{tableprocpower}
\begin{tabular*}{\textwidth}{@{\extracolsep{\fill}}ll@{}}
\hline
\multicolumn{2}{c@{}}{\textbf{Procedures not using the value of} $\bolds{\rho}$} \\
\hline
[Bonf] & the raw Bonferroni procedure \\
{}[LR] & Lehmann Romano's procedure \eqref{criticalvaluesMarkov} \\
{}[AugBonf] & augmentation with $\tau_1=\zeta/m$\\
{}[SimLR] & simultaneous $k$-FWE with \eqref{criticalvaluesMarkov} \\
{}[DimMarkovLR] & diminution with \eqref{equCRS} and $c_{\ell}(x)=x \wt{\tau}_{\ell}$ coming from \eqref
{criticalvaluesMarkov}\\
\hline
\end{tabular*}
\begin{tabular*}{\textwidth}{@{\extracolsep{\fill}}ll@{}}
\hline
\multicolumn{2}{c@{}}{\textbf{Procedures incorporating the value of} $\bolds{\rho}$}\\
\hline
[AugEx] & augmentation with ${\tau}_1$ coming from e.b.d.\\
{}[SimEx] & simultaneous $k$-FWE with $\wt{\tau}_{\ell}$ coming from
e.b.d.\\
{}[Split$1/2$] & new procedure \eqref{newcritval} with $\lambda=1/2$
and $K=2$\\
{}[Split0.95] &new procedure \eqref{newcritval} with $\lambda=0.95$
and $K=2$\\
{}[RWExact] & nonmodified $k$-FWE with $\wt{\tau}_{\ell}$ coming from
e.b.d. \\
{}[DimExEx] & new diminution with \eqref{equCrhoExact} and $c_{\ell}(x)=x \wt{\tau}_{\ell}$ coming from e.b.d.\\
{}[DimGuoLR] & diminution following Theorem~3.8 of \citet{GHS2013} \\
&with $c_{\ell}(x)=x \wt{\tau}_{\ell}$ coming from \eqref
{criticalvaluesMarkov}\\
\hline
\end{tabular*}
\end{table}

As a first illustration, conditions (i), (ii) and (iii) of Theorem~\ref
{thasympfact} are satisfied in the case \eqref{equi-correlated} for
some $\rho\in[0,1)$. Hence, a direct corollary is that the step-up
procedure with the critical values given by \eqref{simpletaul} controls\break 
the FDP asymptotically. Furthermore, in Section S-2 
of the\break \hyperref[suppA]{Supplementary Material}, we complement this result by proving that this FDP
control is maintained if the value of $\rho$ in \eqref{simpletaul} is
replaced by any estimator $\wh{\rho}_m$ provided that
%
%
\begin{equation}
\label{rhochapconsist} (\log m) (\wh{\rho}_m - {\rho})^2
=o_P(1).
\end{equation}

As a second illustration, consider a model \eqref{facmodel} where $c_1$
is uniform on $\{k/r,0\leq k\leq r\}$, and $ \xi_1$ is $\gamma
$-Subbotin (for some $\gamma>1$). Here, $W$ can be any random variable
with continuous distribution function.
For this particular dependence structure, Theorem~\ref{thasympfact}
establishes the asymptotic FDP control of the step-up procedure with
critical values $\tau_{\ell}$ given by the equation
\[
(r+1)^{-1} \sum_{k=0}^{r}
\ol{D}_\gamma \bigl( \ol{F}^{-1}(\tau _{\ell})-
q_\zeta k/r \bigr)= \alpha{\ell}/m,\qquad 1\leq{\ell}\leq m,
\]
where $\ol{D}_\gamma$ denotes the upper-tail function of a $\gamma
$-Subbotin distribution.

\section{Numerical experiments}\label{sec:num}
\label{sec:power}

This section evaluates the power of the procedures considered in
Section~\ref{sec:finitebound} with a proven FDP control. The power is
evaluated by using the standard false nondiscovery rate (FNR), defined
as the expected ratio of errors among the accepted null hypotheses.
Table~\ref{tableprocpower} summarizes the procedures that have been
considered. The simulation are made in model \eqref{equi-correlated}
where the alternative means $\mu_i$ are all equal to some parameter
$\beta$.

%
\begin{figure}

\includegraphics{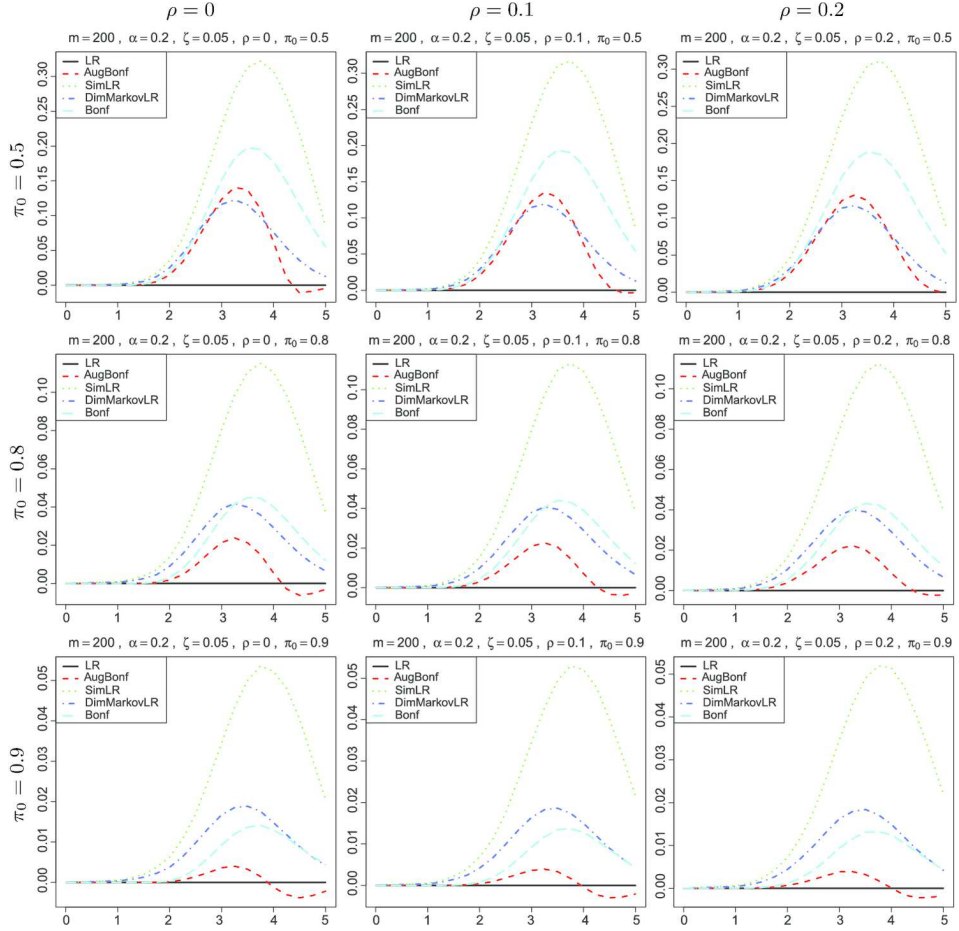}

\caption{Relative FNR to the Lehmann--Romano procedure in function of
$\beta$; see text and Table~\protect\ref{tableprocpower}. Procedures
not using the value of $\rho$.}
\label{powernotusingrho}
\end{figure}

Figure~\ref{powernotusingrho} displays the power of procedures that do
not incorporate the value of~$\rho$; see Table~\ref{tableprocpower}.
Note that, according to Proposition~\ref{RWalreadyproved}(ii),
[LR] controls the FDP because Simes's inequality is valid here.
Hence it does not use the true value of $\rho$ but uses nevertheless an
assumption on the dependence structure. This is not the case of
[AugBonf], [SimLR], [DimMarkovLR] and [Bonf] that
control the FDP for any dependence structure. As one can expect,
[LR] essentially dominates the other procedures. Also, while
[AugBonf] comes in second position, [SimLR] is even worst than
[Bonf] and should be avoided here.

%
\begin{figure}

\includegraphics{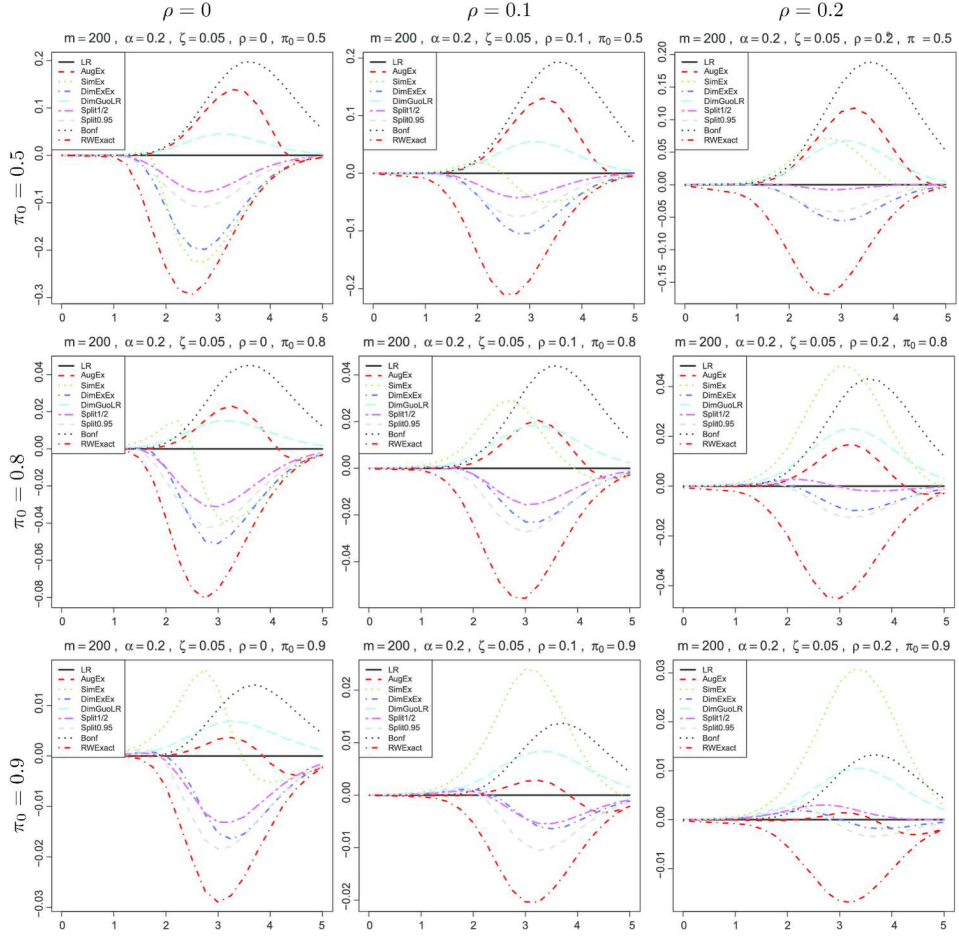}

\caption{Relative FNR to the Lehmann--Romano procedure in function of
$\beta$; see text and Table~\protect\ref{tableprocpower}.
Procedures using the value of $\rho$ (except [LR] and [Bonf]).}
\label{powerusingrho}
\end{figure}

Now, while incorporating the value of $\rho$, we will loosely say that
a procedure is \textit{admissible} if it performs better than [LR] at
least for \textit{a reasonable amount} of parameter configurations.
Figure~\ref{powerusingrho} displays the power of procedures
incorporating the value of $\rho$ (except [LR] and [Bonf] that
we have added only for comparison); see Table~\ref{tableprocpower}.
Note that, except [RWExact], all the procedures have a proven FDP
control, so that the power comparison is fair. First, [DimGuoLR] is
\textit{not admissible}, which indicates that the interest of the bounds
found in \citet{GHS2013} are mainly theoretical in our setting. Second,
[AugEx] only improves [LR] in a very small region, which shows
that, as one can expect, providing $1$-FWE control for controlling the
FDP is too conservative in general. As for [SimEx], things are more
balanced: when $\rho=0$, it improves [LR] when many rejections are
possible ($\pi_0$ not large or $\beta$ large) but does worst otherwise.
We think that this is due to the nature of the [LR] critical
values, which are design to perform well when only few nulls are
expected to be rejected. When $\rho$ is larger, however, [SimEx]
quickly deteriorates. An explanation is that the simultaneity in
[SimEx] is obtained via an union bound, which is conservative when
the dependence is strong. Finally, our new procedures [Split1/2],
[Split0.95] and [DimExEx] seem to be all \textit{admissible} in
these simulations and substantially outperform the other procedures.
Also, none of the three procedures uniformly dominates the others. For
instance, taking $\lambda=1/2$ rather than $\lambda=0.95$ is better
when less rejections are expected, but worst otherwise, while
[DimExEx] seems often better than [Split1/2].

Let us mention that additional simulations have been done in Section
S-5 
of the \hyperref[suppA]{Supplementary Material}. We briefly report some comments here. First, the large
value $\rho=0.5$ has been tried. This deteriorates the relative
performance of all the procedures (except maybe [AugEx]), and in
particular of the $K$-Markov based procedure because the distribution
of the maximum between null $p$-values get closer to the uniform.
Second, simulations have been performed in model \eqref
{alt-equi-correlated} for $a=0.5$. In this model, the positive
dependence property is lost. Hence while [DimExEx] still provably
controls the FDP, this is not anymore the case of [LR],
[Split1/2] and [Split0.95]. However, it is interesting to note
that the FDP control seems to be maintained in the simulations; see
Figure S-4 
in the \hyperref[suppA]{Supplementary Material}. As for the power, the conclusions are qualitatively the
same as in model \eqref{equi-correlated}.

%
%
%
%
%
%
%
%
%
%
%
%
%
%
%
%
%
%
%

%
%
%
%
%
%
%
%
%
%
%
%
%
%
%
%
%
%
%

\section{Conclusion and discussion}\label{sec:discuss}

This paper investigated the FDP control in the case where the
dependence is partly/fully incorporated, by using an extension of RW's
heuristic. We provided two new approaches that offer finite sample
control: the first one (Theorem~\ref{main-th1}) followed the diminution
principle and can be used as soon as the joint distribution of the null
$p$-values can be computed. The second one (Theorem~\ref
{th-newcontrol}) offered a finite sample control under a particular type
of positive dependence \eqref{weakPRD} and exchangeability.
Next, an important part of our work concerned the asymptotic FDP control:
while we established that RW's heuristic is valid under weak dependence
(Theorems~\ref{thRwheuriweakdep} and~\ref{thRwheuriweakdep0}), we
noticed that the interest of the latter has to be balanced with the fact
that the simple BH procedure can be used in this case (Lemma~\ref
{BHFDP}). Then, still based on RW's heuristic, we proposed new critical
values that provide asymptotic control under model \eqref{facmodel}
(Theorems~\ref{thmodifasymp} and~\ref{thasympfact}). Markedly, while it
still relies on a positive dependence assumption \eqref{asympposdep},
this condition has a much simpler form than \eqref{weakPRD}.

Our leading example is related to one-sided testing, so we can
legitimately ask whether our results can be extended to two-sided
testing, that is, when $p_i=2\ol{\Phi}(|X_i|)$ (by using the notation of
Section~\ref{sec:setting}). In the model \eqref{facmodel} with $\xi
_1\sim-\xi_1$, the bounding device calculations done in Section~\ref
{sec:examples} can be clearly generalized to the two-sided case by
replacing $F_{0}(t,w)$ by
$
F_{0}^{(2)}(t,w)= F_{0}(t/2,w)+F_{0}(t/2,-w)$.
Hence we can define new critical values coming from the corresponding
exact bounding device and combine it with the diminution principle
presented in Theorem~\ref{main-th1}. However, the other results of the
paper cannot be directly generalized to the two-sided case because
$F_{0}^{(2)}(t,w)$ may be not increasing w.r.t. $w$.

While this paper solved some issues, it opened several directions of
research. For instance,
is the asymptotic FDP control of Theorem~\ref{thasympfact} still true
when using the original critical values of RW's method rather than
their asymptotic counterpart?
We believe that this issue intrinsically relies on the Poisson
asymptotic regime, which was (essentially) not considered here in our
asymptotic FDP controlling results. Finally, a crucial, but probably
very challenging issue is the validity of RW's approach in the case of
permutation tests
with an arbitrary and unknown dependence structure.

\section{Proofs for finite sample results}\label{sec:proofsfinitesample}

\subsection{An unifying bound} \label{sec:unif}

%
\begin{proposition}\label{prop:finiteboundSUD}
For any critical values $(\tau_{\ell})_{1\leq{\ell}\leq m}$,
consider either
the corresponding step-down \eqref{SD-classic} or step-up \eqref
{SU-classic} procedure, with rejection number $\wh{{\ell}}$. Then the
following holds, both in the fixed model ($\Theta=\Theta^F$) and the
uniform model ($\Theta=\Theta^U$): for all $\theta\in\Theta$,
%
%
\begin{equation}
\P_\theta \bigl(\FDP_m(\tau_{\hat{{\ell}}}) >\alpha \bigr)
\leq\sum_{{\ell}
=1}^{b_\alpha(m_0)} \P_{\theta}
\bigl( V_m(\tau_{{\ell}}) \geq d({\ell},m,m_0),
\wt{{\ell}}={\ell} \bigr),\label{equ-finitebound}
\end{equation}
where $b_\alpha(m_0)$, $d({\ell},m,m_0)$ are defined in \eqref
{balphau} and
\eqref{dlmu}, respectively, and $ \wt{{\ell}}$ is taken as follows:
\begin{longlist}[(ii)]
\item[(i)] Step-down case: $ \wt{{\ell}}=\wh{{\ell}}^{(1)}$, where
$\wh
{{\ell}
}^{(1)}=\min\{{\ell}\in\{1,\ldots,m\} \dvtx S_m(\tau_{\ell})<
(1-\alpha
){\ell}\}$
(with the convention $\min\varnothing=m+1$) and by denoting
$S_m(t)=R_m(t)-V_m(t)$ the number of true discoveries at threshold $t$.
\item[(ii)] Step-up case: $ \wt{{\ell}}=\wh{{\ell}}$.
\end{longlist}
Moreover, in the step-up case, \eqref{equ-finitebound} is an equality.
\end{proposition}

Proposition~\ref{prop:finiteboundSUD}(i) is a reformulation of
Theorem~5.2 in \citet{Roq2011} in our framework and is based on ideas
presented in the proofs of \citet{LR2005,RW2007}.
Proposition~\ref{prop:finiteboundSUD}(ii) is essentially based on
\citet
{RS2006b}, and we provide a short proof below.
\begin{pf*}{Proof of Proposition \ref{prop:finiteboundSUD}}
Since $\FDP_m(\tau_{\hat{{\ell}}}) >\alpha$ implies $\lfloor
\alpha
\wh{{\ell}
}\rfloor+1 \leq m_0$, we have $\wh{{\ell}}\leq b_\alpha(m_0)$.
Also, $ \wh{{\ell}} = R_m(\tau_{\hat{\ell}})\leq m_1+ V_m(\tau
_{\hat{\ell}})$,
which implies $V_m(\tau_{\hat{\ell}}) \geq\wh{{\ell}}-m_1$.
This implies~\eqref{equ-finitebound} in case (ii).
\end{pf*}

\subsection{A new bound}\label{sec:newbound}

%
\begin{proposition}\label{prop:finitebound2}
In the setting of Proposition~\ref{prop:finiteboundSUD}, assume
moreover that there exists a family of random variables $(Z_{{\ell
},{\ell}
'})_{1\leq{\ell},{\ell}'\leq m}$ satisfying: for all ${\ell},{\ell}'$,
%
%
\begin{equation}
\label{equ:Zlk} {\mbf{1}\bigl\{V_m(\tau_{\ell})\geq d\bigl({
\ell}',m,m_0\bigr)\bigr\}}\leq Z_{{\ell
},{\ell}'}\qquad \mbox{a.s.}
\end{equation}
and, a.s., $Z_{{\ell},{\ell}'}$ is nondecreasing in ${\ell}$ and
nonincreasing in
${\ell}'$.
Then for all $\theta\in\Theta$,
%
%
\begin{eqnarray}\label{boundconservSUD}
&&\P_\theta \bigl(\FDP_m(\tau_{\hat{{\ell}}}) >\alpha \bigr)
\nonumber
\\[-8pt]
\\[-8pt]
\nonumber
&&\qquad\leq \sum_{{\ell}
=1}^{b_\alpha(m_0)} \bigl(
\E_\theta ( Z_{{\ell},{\ell}-1} ) -\E_\theta ( Z_{{\ell}-1,{\ell}
-1} )
\bigr) \wedge \bigl(\E_\theta ( Z_{{\ell},{\ell}} ) -\E _\theta
( Z_{{\ell}-1,{\ell}} ) \bigr),
\end{eqnarray}
by letting $Z_{0,{\ell}'}=0$ and $Z_{{\ell},0}=1$ for ${\ell}'\geq
0,{\ell}\geq1$.
\end{proposition}

Applied with $Z_{{\ell},{\ell}'} = V_m(\tau_{\ell})/d({\ell}',m,m_0)$,
Proposition~\ref
{prop:finitebound2} establishes the Romano--Shaikh bound \eqref
{equCRS}. Applied with $Z_{{\ell},{\ell}'} ={\mbf{1}\{V_m(\tau
_{\ell})\geq d({\ell}',m,m_0)\}}$,
Proposition~\ref{prop:finitebound2} entails Theorem~\ref{main-th1}.

\begin{pf*}{Proof of Proposition \ref{prop:finitebound2}}
From \eqref{equ-finitebound}, we derive
\begin{eqnarray*}
\P_\theta \bigl(\FDP_m(\tau_{\hat{{\ell}}}) >\alpha \bigr)
&\leq &\sum_{{\ell}
=1}^{b_\alpha(m_0)} \P_\theta
\bigl( V_m(\tau_{{\ell}}) \geq d({\ell} ,m,m_0) ,
\wt{{\ell}}={\ell} \bigr)\\
&\leq&\sum
_{{\ell}=1}^{b_\alpha(m_0)} \E_\theta \bigl(
Z_{{\ell
},{\ell}} {\mbf{1}\{\wt{{\ell} }={\ell}\}} \bigr).
\end{eqnarray*}
Now, the RHS of the previous display is equal to
\begin{eqnarray*}
&& \sum_{{\ell}=1}^{b_\alpha(m_0)} \E_\theta \bigl(
Z_{{\ell
},{\ell}} {\mbf{1}\{\wt{{\ell} }\geq{\ell}\}} \bigr)-\sum
_{{\ell}=1}^{b_\alpha(m_0)-1} \E_\theta \bigl(
Z_{{\ell}
,{\ell}} {\mbf{1}\{\wt{{\ell}}\geq{\ell}+1 \}} \bigr)
\\
&&\qquad= \sum_{{\ell}=1}^{b_\alpha(m_0)} \E_\theta
\bigl( Z_{{\ell
},{\ell}} {\mbf{1}\{\wt{{\ell} }\geq{\ell}\}} \bigr)-\sum
_{{\ell}=1}^{b_\alpha(m_0)} \E_\theta \bigl(
Z_{{\ell}
-1,{\ell}-1} {\mbf{1}\{\wt{{\ell}}\geq{\ell}\}} \bigr)
\\
&&\qquad= \sum_{{\ell}=1}^{b_\alpha(m_0)} \E_\theta
\bigl( (Z_{{\ell
},{\ell}}- Z_{{\ell}-1,{\ell}
-1}) {\mbf{1}\{\wt{{\ell}}\geq{\ell}\}}
\bigr).
\end{eqnarray*}
Now, since $Z_{{\ell},{\ell}'}$ is nonincreasing w.r.t. ${\ell}'$,
the quantity
$Z_{{\ell},{\ell}}- Z_{{\ell}-1,{\ell}-1}$ is below $ (
Z_{{\ell},{\ell}-1} - Z_{{\ell}
-1,{\ell}-1}
 ) \wedge ( Z_{{\ell},{\ell}} - Z_{{\ell}-1,{\ell}}
 )$, and the latter
is nonnegative because $Z_{{\ell},{\ell}'}$ is nondecreasing w.r.t.
${\ell}$. This
entails the result.
\end{pf*}

\subsection{Proof of Proposition~\texorpdfstring{\protect\ref
{prop:finitebound}}{3.2}}\label
{sec:prooffinitebound}

Let $\wh{k}=V_m(\tau_{\hat{{\ell}}})$, and note that $\wh{k}\leq
m_0$ and
$\{\FDP_m(\tau_{\hat{{\ell}}})>\alpha\}=\{\wh{k}\geq\lfloor
\alpha
\wh{{\ell}}
\rfloor+1\}$.
First, in the nonadaptive case, we have by definition of $\ol{B}_m$,
for all $t$ and $k\leq m_0$,
\[
\ol{B}_m(t ,k )= \sup_{0\leq u \leq m} \bigl\{
B^0_m(t,k,u) \bigr\} \geq B^0_m(t,k,m_0).
\]
Hence, we have by definition of the (nonadaptive) critical values,
\[
\zeta\geq\ol{B}_m\bigl(\tau_{\hat{{\ell}}} , \lfloor\alpha\wh{{
\ell}} \rfloor+1 \bigr) \geq B^0_m\bigl(
\tau_{\hat{{\ell}}},\lfloor\alpha\wh{{\ell}} \rfloor+1,m_0\bigr),
\]
which is larger than or equal to $B^0_m(\tau_{\hat{{\ell}}}, \wh{k},m_0)$
whenever $\wh{k}\geq\lfloor\alpha\wh{{\ell}} \rfloor+1$. Hence
we obtain
\[
\bigl\{\FDP_m(\tau_{\hat{{\ell}}}) >\alpha \bigr\} \subset \bigl\{
B^0_m(\tau_{\hat{{\ell}}}, \wh{k},m_0)\leq
\zeta, \wh{k}\geq1 \bigr\} \subset \bigl\{ \tau_{\hat{{\ell}}} \leq
\nu^0_{\hat{k}} , \wh {k}\geq1 \bigr\},
\]
and thus \eqref{equ-finiteboundSU} holds. Second,
in the adaptive case, we use that $m_0 \leq m - R_m(\tau_{\hat{{\ell
}}}) +
V_m(\tau_{\hat{{\ell}}}) = m-\wh{{\ell}}+\wh{k}$. Thus whenever
$\wh
{k}\geq
\lfloor\alpha\wh{{\ell}} \rfloor+1$, we have for all $t$,
%
%
\begin{eqnarray}\label{equ-proofinterm}
\zeta&\geq&\wt{B}_m\bigl(\tau_{\hat{{\ell}}}, \lfloor\alpha\wh{{\ell}}
\rfloor+1,\wh {{\ell}}\bigr) = \sup_{\lfloor\alpha\hat{{\ell}} \rfloor+1\leq
k'\leq\hat
{{\ell}}} \Bigl\{\sup
_{ 0\leq u \leq m-\hat{{\ell}}+k'} B^0_m\bigl(\tau_{\hat
{{\ell}}},k',u
\bigr) \Bigr\}
\nonumber
\\[-8pt]
\\[-8pt]
\nonumber
& \geq&\sup_{ 0\leq u \leq m-\hat{{\ell}}+\hat{k}} B^0_m(\tau
_{\hat
{{\ell}}},\wh {k},u) \geq B^0_m(
\tau_{\hat{{\ell}}},\wh{k},m_0).
\end{eqnarray}

Hence, this implies $\tau_{\hat{{\ell}}}\leq\nu^0_{\hat{k}}$ and
$\wh
{k}\geq1$, and the proof is complete.

\subsection{Proof of Theorem~\texorpdfstring{\protect\ref{th-newcontrol}}{3.6}}\label
{proofth-newcontrol}

First observe that the critical values \eqref{newcritval} can be
obtained by modifying the $K$-Markov Bounding device $B^0_m(t,k,u)$
defined by \eqref{KMarkovbounding}
as follows:
\[
\bigl(B^0_m\bigr)'(t,k,u)=\cases{
B^0_m(t,k,u)/\lambda,&\quad
$\mbox{if $k\geq K$},$
\vspace*{2pt}\cr
\displaystyle\frac{ut}{(1-\lambda)k}\vee\bigl(B^0_m(t,K,m)/\lambda\bigr),
&\quad
$\mbox{if $k<K$},$}
\]
(the second bounding value being infinite when $\lambda=1$).
Note that the associated adaptive bounding device \eqref
{boundingfunctionadapt} is equal to $(B^0_m)'(t,k,m-{\ell}+k)$ and thus
gives rise to the adaptive critical values \eqref{newcritval}.
By using Proposition~\ref{prop:finitebound} and by letting $\widehat
{k}=V_m(\tau^{\mathrm{new}}_{\hat{{\ell}}})$, we get
\[
\P_\theta \bigl(\FDP_m\bigl(\tau^{\mathrm{new}}_{\hat{{\ell}}}
\bigr) >\alpha \bigr) = \sum_{k=1}^{m_0}
\P_{\theta} \bigl( V_m\bigl(\nu^0_{k}
\bigr) \geq k, \wh{k}=k \bigr),
\]
where $\nu^0_k=\max\{t\in[0,1]\dvtx {B}^0_{m}(t,{k},m_0) \leq
\lambda
\zeta\}$ for all $k\geq K$ and where $\nu^0_k=\max\{t\in
[0,1]\dvtx ({B}^0_{m}(t,{K},m)/\lambda)\vee(m_0 t/(k(1-\lambda))) \leq\lambda
\zeta\}$ for all $k< K$.
It follows that the above display is smaller than or equal to
$T_1+T_2$, where we let
\[
T_1  = \sum_{k=K}^{m_0}
\P_{\theta} \bigl( V_m\bigl(\nu^0_{k}
\bigr) \geq k, \wh {k}=k \bigr) ;\qquad T_2=\sum
_{k=1}^{(K-1)\wedge m_0} \P_{\theta} \bigl( V_m
\bigl(\nu^0_{k}\bigr) \geq k, \wh{k}=k \bigr),
\]
with by convention $T_1=0$ when $K>m_0$.
By \eqref{equ:majKmarkov}, \eqref{exchangeability}, and since $\wh{k}$
is permutation invariant (as a function of the $p$-values), we obtain
\begin{eqnarray*}
T_1 &\leq&\sum_{k=K}^{m_0}
\frac{1}{ {k \choose K} }\sum_{X\subset
{\mtc{H}}
_0\dvtx |X|=K} \P_\theta \Bigl(
\wh{k}=k, \max_{ i\in X}\{ p_i\}\leq\nu
^0_{k} \Bigr)
\\
&=&\sum_{k=K}^{m_0} \frac{m_0(m_0-1)\cdots(m_0-K+1)}{ k(k-1)\cdots
(k-K+1)}
\P_\theta \Bigl(\wh{k}=k, \max_{1\leq i\leq K}\{
q_i\}\leq \nu ^0_{k} \Bigr),
\end{eqnarray*}
where $q_1,\ldots,q_{m_0}$ denotes the $p$-values under the null, that
is, the $p$-values of the set $\{p_i, i\in{\mtc{H}}_0\}$. Next, by
using that
$B^0_m(\nu^0_{k},k,m_0)\leq\lambda\zeta$ for $k\geq K$ and \eqref
{KMarkovbounding}, we get
\begin{eqnarray*}
T_1 &\leq&\lambda\zeta\sum_{k=K}^{m_0}
\P_\theta \Bigl(\wh{k}=k | \max_{1\leq i\leq K}\{
q_i\}\leq\nu^0_{k} \Bigr)
\\
&\leq&\lambda\zeta\sum_{k=K}^{m_0} \Bigl\{
\P_\theta \Bigl(\wh {k}\leq k | \max_{1\leq i\leq K}\{
q_i\}\leq\nu^0_{k} \Bigr) -\P
_\theta \Bigl(\wh{k}\leq k-1 | \max_{1\leq i\leq K}\{
q_i\}\leq\nu ^0_{k} \Bigr) \Bigr\}.
\end{eqnarray*}
Now, since the $p$-value subset of $[0,1]^m$ defined by the relation
$\wh{k}\leq k-1$ is nondecreasing, assumption~\eqref{weakPRD} ensures
\begin{eqnarray*}
T_1 \leq\lambda\zeta\sum_{k=K}^{m_0}
\Bigl\{ \P_\theta \Bigl(\wh {k}\leq k | \max_{1\leq i\leq K}\{
q_i\}\leq\nu^0_{k} \Bigr) -\P
_\theta \Bigl(\wh{k}\leq k-1 | \max_{1\leq i\leq K}\{
q_i\} \leq\nu ^0_{k-1} \Bigr) \Bigr\} ,
\end{eqnarray*}
which is below $ \lambda\zeta$ because the sum is telescopic.

Now, for $T_2$, we use the same type of reasoning with ${\mbf{1}\{
V_m(\nu ^0_{k}) \geq k\}}\leq\frac{1}{ k} \sum_{i=1}^{m_0} {\mbf
{1}\{q_i\leq\nu ^0_{k}\}}$ and ${m_0 \nu^0_{k}}\leq(1-\lambda
)\zeta k$ for $k<K$,
\begin{eqnarray*}
T_2 &\leq&\sum_{i=1}^{m_0} \sum
_{k=1}^{(K-1)\wedge m_0} \frac{1}{ k} \P
_\theta \bigl(\wh{k}=k, q_i\leq\nu^0_{k}
\bigr)
\\
&\leq&(1-\lambda) \zeta m_0^{-1}\sum
_{i=1}^{m_0}\sum_{k=1}^{(K-1)\wedge m_0}
\P_\theta \bigl(\wh{k}=k | q_i\leq\nu ^0_{k}
\bigr)
\\
&\leq&(1-\lambda) \zeta m_0^{-1}\sum
_{i=1}^{m_0}\sum_{k=1}^{(K-1)\wedge m_0}
\P_\theta \bigl(\wh{k}\leq k | q_i\leq \nu
^0_{k} \bigr) -\P_\theta \bigl(\wh{k}\leq k-1 |
q_i\leq\nu ^0_{k} \bigr)
\\
&\leq&(1-\lambda) \zeta m_0^{-1}\sum
_{i=1}^{m_0}\sum_{k=1}^{(K-1)\wedge m_0}
\bigl\{ \P_\theta \bigl(\wh{k}\leq k | q_i\leq
\nu^0_{k} \bigr) -\P_\theta \bigl(\wh{k}\leq k-1
| q_i\leq \nu^0_{k-1} \bigr) \bigr\},
\end{eqnarray*}
by using again assumption~\eqref{weakPRD}. Finally, the last display is
below $(1-\lambda) \zeta$, because the sum is telescopic.
This completes the proof.

\section{Proofs for asymptotic results}\label{sec:proofs}

In this section, the following well-known lemma will be extensively used.

%
\begin{lemma}\label{lemmenullquantileproc}
Let $\wh{{\ell}}$ be the number of rejections of the step-up \eqref
{SU-classic} algorithm associated to
some critical values $(\tau_{\ell})_{1\leq{\ell}\leq m}$.
Consider the function $f_m$ defined by
%
%
\begin{equation}
f_m(t)=m^{-1}\times\min\bigl\{{\ell}\in\{0,\ldots,m+1\}\dvtx \tau_{\ell}
 \geq t\bigr\}, \label{fmt}
\end{equation}
with the conventions $\tau_0=0$, $\tau_{m+1}=1$. Let $\wh{t}$ be
defined by
%
%
\begin{equation}
\label{tchap} \wh{t} = \sup\bigl\{t\in[0,1] \dvtx{\widehat{
\mathbb{G}}}_{m}(t)\geq f_m(t)\bigr\}.
\end{equation}
Then the supremum into \eqref{tchap} is a maximum, that is,
${\widehat{\mathbb{G}}}_{m}(\wh{t}) \geq f_m(\wh{t})$. Furthermore,
$\wh{t}=\tau_{\hat
{{\ell}}}$.
\end{lemma}

\subsection{Proof of Lemma~\texorpdfstring{\protect\ref{BHFDP}}{4.1}}\label{sec:proofBHFDP}

Actually, we prove the result for a more general class of procedures,
where $\wh{t}=\tau_{\hat{\ell}}$ is obtained by \eqref{tchap} for a
sequence of functions $f_m=\wh{f}_m$ (possibly random) which is
uniformly close to $f_\infty(t)=t/\alpha$ on every compact of
$(0,\alpha
]$, that is,
%
%
\begin{equation}
\label{unif} \sup_{b\leq t\leq\alpha}\bigl|\wh{f}_m(t)-t/\alpha\bigr|
\rightarrow0 \qquad\mbox{a.s. for all $b\in(0,\alpha)$.}
\end{equation}
Note that
$\wh{f}_m=f_\infty$ gives the BH procedure by Lemma~\ref
{lemmenullquantileproc}.
Next, since $R_m(\wh{t})\geq m \wh{f}_m(\wh{t})$,
\begin{eqnarray*}
&&\P \bigl( \FDP_m (\wh{t} )> \alpha \bigr)\\
&&\qquad\leq \P \bigl(
V_m(\wh {t})/m > \alpha\wh{f}_m(\wh{t}) \bigr)
\\
&&\qquad= \P \bigl( (m_0/m) \bigl({\widehat{\mathbb{G}}}_{0,m}(
\wh{t}) -\wh{t} \bigr) - \alpha\bigl(\wh {f}_m(\wh {t})- \wh{t}/\alpha
\bigr) > (1-m_0/m) \wh{t} \bigr)
\\
&&\qquad\leq \P \Bigl( (m_0/m) \Vert{\widehat{\mathbb{G}}}_{0,m}-
I \Vert _\infty+ \alpha\sup_{t^\star
\leq t\leq\alpha}\bigl|
\wh{f}_m(t)-t/\alpha\bigr|> (1-m_0/m) t^\star \Bigr)
\\
&&\qquad\quad{} + \P\bigl(\wh{t}\leq t^\star\bigr),
\end{eqnarray*}
for some $t^\star>0$ satisfying $G(t^\star)>t^\star/\alpha$ [which
exists by \eqref{exists}]. By \eqref{condm0} and \eqref{weakdep0}, it
is sufficient to check that $\P(\wh{t}\leq t^\star)$ tends to zero. For
this, we use \eqref{weakdep} that ensures
\begin{eqnarray*}
\P\bigl(\wh{t}> t^\star\bigr) &\geq&\P\bigl({\widehat{
\mathbb{G}}}_{m}\bigl(t^\star \bigr)> \wh{f}_m
\bigl(t^\star\bigr) \bigr)
\\
&\geq&\P\bigl(G\bigl(t^\star\bigr)>{f}_\infty
\bigl(t^\star\bigr) +\bigl| G\bigl(t^\star\bigr)- {\widehat {
\mathbb{G}}}_{m} \bigl(t^\star\bigr)\bigr| + \bigl|\wh{f}_m
\bigl(t^\star\bigr) - {f}_\infty\bigl(t^\star\bigr)\bigr|
\bigr) \rightarrow1,
\end{eqnarray*}
which completes the proof.

\subsection{Proof of Theorem~\texorpdfstring{\protect\ref{thRwheuriweakdep0}}{4.2}}\label
{proof:thRwheuriweakdep0}

By the proof of Lemma~\ref{BHFDP}, it is sufficient to show that
$f_m(t)$ defined by \eqref{fmt} is such that $f_m(t)\rightarrow
t/\alpha
$ for all $t\in[0,1]$.
This is an easy consequence of the fact that, since $V'_m(t)/m$
converges in probability to $t$, for any sequence $({\ell}_m)_m$ with
${\ell}
_m/m$ converging to some $u$, $\ol{B}_m(t,\lfloor\alpha{\ell}_m
\rfloor
+1)$ converges to~$1$ if $ t>\alpha u$ and $0$ if $ t<\alpha u$.

\subsection{Proof of Theorem~\texorpdfstring{\protect\ref{thRwheuriweakdep}}{4.3}}\label
{sec:thRwheuriweakdep}

First, by assumption~\eqref{weakdep}, we can assume that the
convergence $\sup_{t\in[0,1]}|{\widehat{\mathbb
{G}}}_{m}(t)-G(t)|\rightarrow0$ is
almost sure.
Next, let us prove
%
%
\begin{equation}
\label{convtstar} \wh{t} \mbox{ converges a.s. to } t^\star\in(0,1),
\end{equation}
where ${t}^{\star}=\sup\{t\in[0,1]\dvtx G(t)\geq\pi_0 t/\alpha\}$.
First, $t^\star$ lies in $(0,1)$ by \eqref{noncritical} and because
$\pi
_0>\alpha$.
Then, by Lemma~\ref{lemmenullquantileproc}, we have $\wh{t} = \sup\{
t\in
[0,1] \dvtx{\widehat{\mathbb{G}}}_{m}(t)\geq f_m(t)\}$ where
$f_m(t)$ is given \eqref{fmt}.
As in proof of Theorem~\ref{thRwheuriweakdep0}, we easily check that
for all $t\in[0,1]$, $f_m(t)$ converges to $\pi_0 t/\alpha$.
As a result, since $f_m$ is a nondecreasing function, the convergence
of $f_m(t)$ to $\pi_0 t/\alpha$ is uniform on $[0,1]$.
Now, to establish~\eqref{convtstar}, it is sufficient to show that if
$\wh{t}$ converges to some $t\in[0,1]$ along a subsequence, then we
have $t={t}^{\star}$.
First, since ${\widehat{\mathbb{G}}}_{m}(\wh{t})\geq f_m(\wh{t})$,
we have $G(t)\geq\pi_0
t/\alpha$ and thus $t\leq t^\star$. Let us prove $t\geq t^\star$.
We have by \eqref{unique} and \eqref{noncritical} that $G(u_p)>
f_\infty
(u_p)$ for all $p$,\vspace*{1pt} for some $u_p \uparrow t^\star$. This yields, for
all $p$ and $m$ large enough, ${\widehat{\mathbb{G}}}_{m}(u_p)>
f_m(u_p)$ and thus $t\geq
u_p$. Hence, $t\geq t^\star$ by making $p$ tends to infinity. This
proves \eqref{convtstar}.

Now, we have
$
\P (\FDP_m(\wh{t}) >\alpha ) = \P (V_m(\wh{t})
>\alpha\wh
{{\ell}} ) = \P (Z_m(\wh{t}) > \Upsilon_m )$,
by letting $\Upsilon_m=r_m( \alpha{\wh{{\ell}}}/m - \tau^0_{\hat
{{\ell}}}
m_0/m )$. By assumption~\eqref{FCLT}, we have that $Z_m(\wh{t})$
converges in distribution to $Z({t}^\star)$. Let $q_m(t)$ denotes the
$(1-\zeta)$-quantile of $Z_m(t)$.
From Lemma S-3.2, 
we have that the function sequence $q_m(t)$ converges uniformly to
$q_\zeta(t)$ for $t$ in any compact of $(0,1)$, where $q_\zeta(t)$
denotes the $(1-\zeta)$-quantile of $Z(t)$.
From above, the proof is complete if we show
%
%
\begin{equation}
\label{convUpsi} \Upsilon_m \mbox{ converges a.s. to }
q_\zeta\bigl(t^\star\bigr).
\end{equation}

Let us prove \eqref{convUpsi}. By definition of $B_m^0$, we have $\P
(V_m(\tau^0_{\ell})>\alpha{\ell}) \leq\zeta< \P(V_m((\tau
^0_{\ell}
+\varepsilon
/r_m)\wedge1)>\alpha{\ell})$, for all $\varepsilon>0$.
Note that the latter uses that $m_0>\alpha m$ (for $m$ large enough).
This shows that for all ${\ell}\in\{1,\ldots,m\}$,
$q_m(\tau^0_{\ell}) \leq r_m( \alpha{{\ell}}/m - \tau^0_{\ell
}m_0/m )\leq
q_m((\tau^0_{\ell}+\varepsilon/r_m)\wedge1)+\varepsilon$. 
Hence, applying this relation to ${\ell}=\wh{{\ell}}$, we get that
for all
$\varepsilon>0$, a.s.,
$q_\zeta(t^\star)\leq\liminf_m \Upsilon_m\leq\limsup_m \Upsilon
_m\leq
q_\zeta(t^\star)+\varepsilon$. Then~\eqref{convUpsi} is derived by
making $\varepsilon$ tend to zero.

\subsection{Proof of Theorem~\texorpdfstring{\protect\ref{thmodifasymp}}{4.4}}\label
{sec:proofwithDKW}

We have
%
%
\begin{eqnarray}\label{equintermtheorem44}
\P \bigl(\FDP_m(\wh{t}) >\alpha \bigr)& \leq& \P
\bigl((m_0/m) {\widehat{\mathbb{G}}}_{0,m}(\wh {t}) >\alpha
\wh{{\ell}}/m ,\wh{t}>0 \bigr)
\nonumber
\\
&\leq& \P \biggl( \bigl\Vert{\widehat{\mathbb{G}}}_{0,m}(\cdot) -
F_0(\cdot,W) \bigr\Vert_\infty> \biggl(\frac
{-\log(\lambda\zeta/2)}{2m_0}
\biggr)^{1/2} \biggr)
\\
&&{}+ \P \bigl( F_0(\wh{t},W) > F_{0}(\wh{t},q_{\zeta(1-\lambda)})
,\wh {t}\in(0,1) \bigr) .\nonumber
\end{eqnarray}
Now, conditionally on $W$, the $p_i$'s are i.i.d. of distribution
function $F_0(\cdot,W)$. Hence, by applying the
Dvoretzky--Kiefer--Wolfowitz inequality with the tight constant [see
\citet{Mass1990}], we get that the first term in the previous display is
smaller than
$ \lambda\zeta$,
which in turn implies that \eqref{equintermtheorem44} is smaller than
$\lambda\zeta+\P ( W\geq q_{\zeta(1-\lambda)}  ) \leq
\zeta$
by \eqref{asympposdep}.

\subsection{Proof of Theorem~\texorpdfstring{\protect\ref{thasympfact}}{4.5}}\label
{sec:proofthasympfact}

First note that since $\ol{F}_{\xi}$ is continuous and increasing, so
is $\ol{F}$ and $t\in[0,1]\mapsto F_0(t,w)$, for all $w$. Hence
\eqref
{equ-critvalues-RWasymp} defines the $\tau_{\ell}$'s in an unique manner.
Next, if $\P(c_1>0)=0$, then $F_0(t,w)=t$ for all $w$, and thus the
considered procedure is the BH procedure, which controls the FDP
asymptotically by Lemma~\ref{BHFDP}. Hence we can assume that $\P
(c_1>0)>0$. Let us denote the support of $c_1$ by $\{v_1,\ldots,v_r\}$,
for $r\geq1$, $v_i\geq0$, $v_i\neq v_j$ for $i\neq j $. We thus have
that at least one $v_i$ is positive. In particular, assumption~\eqref
{asympposdep} holds.

Then, by using the Skorokhod representation theorem, up to consider a
subsequence, we can assume that $(\wh{t},W)$ is almost surely
converging to some $(T,W)$ (on appropriate subspaces).
Denote $\kappa_\zeta=\break \max_{1\leq i\leq r}\{v_i q_\zeta\}$, $\kappa
_W=\max_{1\leq i\leq r}\{v_i W\}$, and let us establish
%
%
\begin{equation}
\label{bornesup}  T>0\qquad\mbox{a.s. if } \kappa_W+\beta>
\kappa_\zeta.
\end{equation}
For this, note that by Lemma~\ref{lemmenullquantileproc}, $\wh{t}$ is
obtained by \eqref{tchap} with $f_m(t)=F_{0}(t ,q_\zeta)/\alpha$,
which gives
$F_{0}(\wh{t},q_\zeta) = \max\{ t'\in[0,1] | {\widehat{\mathbb
{G}}}_{m}'(t')\geq
t'/\alpha\}
$,
where\break 
${\widehat{\mathbb{G}}}_{m}'(t') =  m^{-1} \sum_{i=1}^m {\mbf{1}\{
F_{0}(p_i,q_\zeta) \leq t'\}} $.
Now observe that there exists a constant $D\in(0,1)$ such that for all $u$,
%
%
\begin{equation}
\label{encadrF0} D \ol{F}_\xi\bigl(\ol{F}^{-1}(u)-
\kappa_\zeta\bigr) \leq F_{0}(u,q_\zeta )\leq\ol
{F}_\xi\bigl(\ol{F}^{-1}(u)-\kappa_\zeta\bigr).
\end{equation}
It follows that ${\widehat{\mathbb{G}}}_{m}'(t') $ is lower-bounded
by $m^{-1} \sum_{i=1}^m
H_i \mbf{1}\{\ol{F}_\xi(\ol{F}^{-1}(p_i)-\kappa_\zeta) \leq  t'\}$,
which by
the law of large numbers [because $(c_i,\xi_i)$ are i.i.d.] converges
a.s. toward
\[
\pi_1 \P \bigl(\ol{F}_\xi\bigl(\ol{F}^{-1}_\xi
\bigl(t'\bigr)+\kappa_\zeta -c_1W-\beta \bigr)
\cond W \bigr)\geq D' \ol{F}_\xi\bigl(
\ol{F}^{-1}_\xi\bigl(t'\bigr)+\kappa
_\zeta- \kappa_W-\beta\bigr),
\]
where $D'$ is some positive constant.
Assume now $\kappa_W+\beta>\kappa_\zeta$. By \eqref{LRseparate}, the
slope of $\ol{F}_\xi(\ol{F}^{-1}_\xi(t')+\kappa_\zeta- \kappa
_W-\beta)$
is infinite in $0$.
Hence, for $m$ large enough we have $F_{0}(\wh{t},q_\zeta)>t'_0$,
where $t'_0$ denotes any $t'\in(0,1)$ such that $D' \ol{F}_\xi(\ol
{F}^{-1}_\xi(t')+\kappa_\zeta- \kappa_W-\beta)>t'/\alpha$. As a result,
$T>0$ and \eqref{bornesup} is proved.

Now, we establish
%
%
\begin{eqnarray}
\label{borneinf} \mbox{For all }
\varepsilon>0, \exists t_\varepsilon(W)
\in(0,1), \mbox{ s.t.}\qquad \P\bigl(\wh{{\ell}}\geq1,
\wh{t}\leq t_\varepsilon(W)\cond W
\bigr)\leq \varepsilon\alpha/D
\nonumber
\\[-8pt]
\\[-8pt]
\eqntext{\mbox{if }\kappa_W+\beta<
\kappa_\zeta.}
\end{eqnarray}
By the LHS of \eqref{encadrF0}, we obtain that ${\widehat{\mathbb
{G}}}_{m}'(t') $ is
upper-bounded by
\begin{eqnarray*}
&&m^{-1} \sum_{i=1}^m {\mbf{1}
\bigl\{\ol{F}_\xi(\beta+c_iW+\xi _i-\kappa
_\zeta) \leq t'/D\bigr\}}
\\
&&\qquad\leq m^{-1} \sum_{i=1}^m {
\mbf{1}\bigl\{q_i\leq\ol{F}_\xi\bigl(\ol
{F}^{-1}_\xi \bigl(t'/D\bigr)+
\kappa_\zeta-\kappa_W-\beta\bigr)\bigr\}},
\end{eqnarray*}
where we let $q_i=\ol{F}_\xi(\xi_i)$, for $1\leq i\leq m$, which are
i.i.d. uniform. Now assume $\kappa_\zeta-\kappa_W-\beta>0$, and take
any $\varepsilon>0$. By \eqref{LRseparate}, there exists
$t'_\varepsilon
(W)\in(0,1)$ such that $\forall t'\in(0,t'_\varepsilon(W)]$, we have
$\ol{F}_\xi(\ol{F}^{-1}_\xi(t'/D)+\kappa_\zeta-\kappa_W-\beta
)\leq
\varepsilon t'/D$.
Then we have
\[
\P\bigl(\wh{{\ell}}\geq1, \alpha\wh{{\ell}}/m\leq t'_\varepsilon
(W)\cond W\bigr) \leq\P \bigl(\exists{\ell}\in\{1,\ldots,m\} \dvtx
q_{({\ell})} \leq \varepsilon D^{-1}\alpha{\ell}/m \cond W \bigr),
\]
which is below $\varepsilon D^{-1}\alpha$ by using Simes's inequality;
see, for example, \eqref{equsimes}. This provides \eqref{borneinf} by
taking $t_\varepsilon(W) \in(0,1)$ such that $F_0(t_\varepsilon
(W),q_\zeta)=t'_\varepsilon(W)$.

The last argument is that when $T>0$ a.s., we have $F_{0}(T,q_\zeta)>0$
a.s. and thus
%
%
\begin{eqnarray}
\FDP(\wh{t}) = \frac{m_0}{m}\frac{{\widehat{\mathbb
{G}}}_{0,m}(\wh{t})}{{\widehat{\mathbb{G}}}_{m}(\wh
{t})} = \frac{m_0}{m}\alpha
\frac{{\widehat{\mathbb{G}}}_{0,m}(\wh
{t})}{F_{0}(\wh{t},q_\zeta)} \rightarrow\pi_0 \alpha\frac{F_{0}(T,W)}{F_{0}(T,q_\zeta)}.\label
{lastargument}
\end{eqnarray}

Now, by combining \eqref{bornesup}, \eqref{borneinf} and \eqref
{lastargument}, we obtain
\begin{eqnarray*}
&&\limsup_m \P\bigl(\FDP_m(\wh{t}) >\alpha\bigr)
\\
&&\qquad\leq \E \Bigl( \limsup_m {\mbf{1}\bigl\{\FDP _m(
\wh{t}) >\alpha, \kappa_W+\beta>\kappa_\zeta\bigr\} }
\Bigr)
\\
&&\qquad\quad{} + \limsup_m \P \bigl(\FDP_m(\wh{t}) >\alpha,
\kappa_W+\beta <\kappa _\zeta, \wh{t}\leq
t_\varepsilon(W) \bigr)
\\
&&\qquad\quad{}+ \E \Bigl(\limsup_m {\mbf{1}\bigl\{\FDP_m(
\wh{t}) >\alpha, \kappa _W+\beta<\kappa _\zeta, \wh{t}>
t_\varepsilon(W) \bigr\}} \Bigr)
\\
&&\qquad\leq \P \biggl(\pi_0 \alpha\frac{F_{0}(T,W)}{F_{0}(T,q_\zeta)} \geq \alpha,
\kappa_W+\beta>\kappa_\zeta, T>0 \biggr) + \varepsilon
\alpha D^{-1}
\\
&&\qquad\quad{}+\P \biggl(\pi_0 \alpha\frac{F_{0}(T,W)}{F_{0}(T,q_\zeta)} \geq \alpha ,
\kappa_W+\beta<\kappa_\zeta, T\geq t_\varepsilon(W)
\biggr).
\end{eqnarray*}
Also note that $T<1$ a.s. on the two above events, because $\pi
_0\alpha
<\alpha$.
Hence we get
\[
\limsup_m \P\bigl(\FDP_m(\wh{t}) >\alpha\bigr)
\leq\P \bigl( F_{0}(T,W) > F_{0}(T,q_\zeta), T
\in(0,1) \bigr) + \varepsilon\alpha D^{-1},
\]
and the result comes from \eqref{asympposdep} and by letting
$\varepsilon$ tends to zero.

\section*{Acknowledgments}
We acknowledge an Associate Editor and three referees for their
relevant suggestions that have substantially improved the quality of
the manuscript.
We are very grateful to Yoav Benjamini, Gilles Blanchard and Pierre
Neuvial for helpful discussions. We also thank Bertrand Michel and
Jacques Portes for their support during the numerical experiments.


\begin{supplement}[id=suppA]
\stitle{Supplement to ``New procedures controlling the false discovery
proportion via Romano--Wolf's heuristic''}
\slink[doi]{10.1214/14-AOS1302SUPP} 
\sdatatype{.pdf}
\sfilename{aos1302\_supp.pdf}
\sdescription{The supplement presents additional materials for the paper; see \citet
{DR2014supp}.}
\end{supplement}

%
%

%

\printaddresses
\end{document}